# Width distributions for rectangular boxes

*A global density, singularity structure and moment closure*

Omri Abas

*Independent researcher, Kadima-Zoran, Israel*

omri.abas1@gmail.com

ORCID: 0009-0003-4857-0790

---

## Abstract

For a rectangular box with edges $a_1, a_2, a_3$, the width in a uniform random direction $u \in S^2$ is $w = \sum_i a_i\,|u_i|$. Under a bijective change of parameters this is also the projected area of a rectangular parallelepiped, whose distribution was derived by Walters. We give a direct co-area derivation that recasts the law as one global positive-part formula and integrates it to an elementary closed cumulative distribution function.

The density is real-analytic except at the edge lengths, face diagonals and space diagonal. We classify every singularity as a corner, square-root fold, superposition or terminal jump and compute its coefficient. A Gaussian representation gives all moments from one generating function. The degree of the $n$th cumulant in $\pi^{-1}$ is exactly $\lfloor n/2 \rfloor$; at every order the mean-normalised cumulant is a polynomial in two scale-invariant combinations of the intrinsic volumes, with explicit formulas given through fifth order. The first three cumulants determine the box, and we characterise the admissible triples.

Beyond boxes, intrinsic volumes determine neither width variance nor brightness variance. The sign-pattern representation extends to bodies whose central symmetral is a zonotope, while the cellwise co-area method applies to every full-dimensional polytope.

---

# 1. Introduction

## 1.1 The width as a random variable

Let $K \subset \mathbb{R}^3$ be a convex body with support function $h_K(u) = \max_{x \in K} x \cdot u$. The **width** of $K$ in a direction $u \in S^2$ is the distance between its two supporting planes normal to $u$,

$$w_K(u) = h_K(u) + h_K(-u)$$

Taking $u$ uniform on $S^2$ makes $w_K$ a random variable. Its expectation is the classical mean width, which together with volume and surface area exhausts the non-trivial intrinsic volumes of a body in $\mathbb{R}^3$, and which has been computed for many bodies since Cauchy. Its **distribution** has received far less attention.

## 1.2 A question of Finch

Finch [4] computes, for the cube of unit edge,

$$\mathbb{E}[w_{\square}] = 3/2, \qquad \mathbb{E}[w_{\square}^2] = 1 + 4/\pi$$

together with the corresponding values for the regular tetrahedron and for the square and the equilateral triangle in the plane, describing the second moments as apparently new. He observes that *the probability density of* $w$ *is not known*, and closes by hoping that the computation might be carried out for other convex bodies. The same random variable had, however, appeared earlier in the projected-area literature.

Walters [16] derived the projected-area distribution of a randomly oriented rectangular parallelepiped. If its edge lengths are $b_1, b_2, b_3$, then its projected area normal to $u$ is

$$A(u) = b_2 b_3 |u_1| + b_3 b_1 |u_2| + b_1 b_2 |u_3|.$$

Conversely, for arbitrary $a_1, a_2, a_3 > 0$, set

$$b_1 = \sqrt{\frac{a_2 a_3}{a_1}}, \qquad b_2 = \sqrt{\frac{a_1 a_3}{a_2}}, \qquad b_3 = \sqrt{\frac{a_1 a_2}{a_3}}.$$

Then $b_2 b_3 = a_1$, $b_3 b_1 = a_2$ and $b_1 b_2 = a_3$, so the projected-area law of rectangular parallelepipeds is exactly the three-parameter law of $\sum_i a_i |u_i|$. Vickers [14] later expressed the distribution for a general cylinder through the width function of its base in a form suitable for numerical computation and gave analytic results for the cube. Umhauer and Gutsch [13] subsequently calculated and measured projected-area frequency distributions for a cube and a rectangular parallelepiped, while Vickers and Brown [15] treated further families of convex particles.

Walters's density is casewise. Writing the three face areas as $A_x \le A_y \le A_z$, he treats $A_x^2 + A_y^2 < A_z^2$ and $A_x^2 + A_y^2 > A_z^2$ separately and divides each regime at the face areas and the pairwise square roots, obtaining inverse-sine formulas on six intervals [16, pp. 345–346]. His $G(A)$, termed a frequency distribution function, is the density; he also records the first three raw moments and the variance. In the present ordering $(A_x, A_y, A_z) = (a_3, a_2, a_1)$, his separating value $\sqrt{A_x^2 + A_y^2}$ is $d_1$, so the two case systems are exactly $d_1 < a_1$ and $d_1 > a_1$; their six intervals are cut by $a_3, a_2, \min(a_1, d_1), \max(a_1, d_1), d_2, d_3, d$. Theorem 3.1 encodes both regimes, every ordering and all coincidences in one symmetric positive-part formula, while Proposition 3.11 supplies an elementary cumulative distribution function. The global form then supports the complete singularity classification, the all-order moment and cumulant algebra, inverse reconstruction and the extensions through central symmetrals developed below.

The problem also arose independently in room acoustics, through the study of directional geometric statistics for rectangular rooms: the full width distribution retains information that the mean width does not. The present paper isolates the underlying geometric question; the applications are developed separately.

Finch's conjectures on the projection length of the regular cube, simplex and crosspolytope were subsequently established by Kabluchko, Litvak and Zaporozhets [7], using a general Gaussian width-moment identity valid for every compact convex body and applying it to regular polytopes in every dimension. Lemma 4.1 below is its direct three-dimensional rectangular-box specialisation; compare [7, Theorem 3.1]. The present paper uses that device for the all-order moment algebra, cumulant structure and inverse reconstruction of the three-parameter family.

Equivalently, by the substitution $x_i = u_i^2$, which pushes the uniform measure on $S^2$ forward to the Dirichlet distribution with parameters $(1/2\,,1/2\,,1/2)$, the random variable studied here is

$$w = \sum_{i=1}^{3} a_i\,\sqrt{x_i}, \qquad x \sim \mathrm{Dir}(1/2\,,1/2\,,1/2)$$

a weighted sum of square roots of a Dirichlet vector. The global formula below and its elementary antiderivative provide the analytic starting point for the subsequent singularity, moment and reconstruction theory.

## 1.3 One measure, two weightings

Every convex body carries its surface area measure $S_K$ on $S^2$; for a polytope with faces of area $S_k$ and outer normals $n_k$ this is the atomic measure $S_K = \sum_k S_k\,\delta_{n_k}$. Two classical quantities are integrals against it, differing only in how the normal is weighted:

$$V(K) = \frac{1}{3}\int_{S^2} h_K\; dS_K, \qquad\qquad A_K(u) = \frac{1}{2}\int_{S^2} |\,u\cdot n|\; dS_K(n)$$

The first is the decomposition of $K$ into pyramids over its faces; the second is Cauchy's projection formula [3], giving the area of the shadow of $K$ normal to $u$; see [12, §5.3]. The distinction is in what is produced rather than in what is integrated: the first has an integrand depending on the variable of integration alone and returns a number, independently of the choice of apex by Minkowski's relation [9], whereas the second has a kernel $|u\cdot n|$ in two variables and returns a function on the sphere. Randomising the free variable of that function is the subject of what follows.

The second integral is the cosine transform of $S_K$, and $A_K$ is by definition the support function of the projection body $\Pi K$ [5, §4.1], [12, §10.9]. Since $\Pi K$ is centrally symmetric,

$$A_K(u) = h_{\Pi K}(u) = 1/2\; w_{\Pi K}(u)$$

for **every** convex body: brightness distributions and width distributions are linked through the projection-body operator. What distinguishes the box is not this identity but that its class is stable under $\Pi$. For $K$ the box with edges $a_1, a_2, a_3$,

$$A_K(u) = a_2 a_3 |u_1| + a_3 a_1 |u_2| + a_1 a_2 |u_3| = V\sum_i \frac{|u_i|}{a_i}$$

so that

$$\Pi K = 2V \cdot K^{\#}, \qquad K^{\#} \text{ the box with edges } 1/a_i \tag{1}$$

The projection body of a box is a box. The notation $K^{\#}$ is used for the reciprocal-edge box to avoid the conventional meaning of $K^*$ as a polar body; the polar of a box is a cross-polytope, not a box. Consequently the width and brightness distributions of boxes form a single family, closed under $a \mapsto 1/a$, and in either reading the random variable is a **linear form in** $|u|$. Equation (1) is the projection-body form of the parameter correspondence displayed in Section 1.2. Every closed form in this paper descends from that linearity: the level sets on $S^2$ are circles, the constraints defining an orthant are three great circles, and the moments are symmetric functions of the edges.

The ball and the ellipsoid are $\Pi$-stable as well. For an ellipsoid with semiaxes $\alpha_i$ the width is $2(\sum_i \alpha_i^2 u_i^2)^{1/2}$, the square root of a **quadratic** form rather than a linear form in $|u|$, and its distribution on the sphere is given piecewise by Hillier [6] on the intervals between the characteristic roots. The comparison is taken up in Question 8.7.

## 1.4 Results

Throughout, $K$ is the box with edge lengths $a_1 \ge a_2 \ge a_3 > 0$ and $w = w_K$ its width at a direction $u$ uniform on $S^2$. We write $d = (\sum_i a_i^2)^{1/2}$ for the space diagonal and $d_i = (d^2 - a_i^2)^{1/2}$ for the diagonal of the face normal to $e_i$. By Remark 2.4 the width takes values in $[a_3, d]$, from the shortest edge to the diameter.

**Theorem A** (Theorem 3.1). *For $\rho \in (a_3, d)$ put $s = \sqrt{d^2 - \rho^2}$ and*

$$t_i = \frac{\rho\, a_i}{s\, d_i}, \qquad \eta_i = \arccos\big(\min(t_i, 1)\big), \qquad \Theta_{ij} = \arccos\left(\frac{-a_i a_j}{d_i d_j}\right)$$

*Then $w$ has density*

$$f(\rho) = \frac{2}{\pi d}\left[2\pi - 2\sum_{i=1}^{3} \eta_i + \sum_{i<j} (\eta_i + \eta_j - \Theta_{ij})_+\right]$$

*on that interval, and $f = 0$ outside* $[a_3, d]$.

Under the parameter correspondence of Section 1.2, Walters [16] gives the same law in a casewise projected-area form. Theorem A is the global width-coordinate form used throughout the present paper: its positive parts encode every ordering of the edge and face-diagonal values in one expression. Proposition 3.11 integrates it to a closed CDF.

**Theorem A′** (Proposition 3.11). *There is an explicit elementary function $\mathcal{P}$, formed from algebraic functions, arccosines and arctangents, such that $\mathcal{P}(a_3) = 0$, $\mathcal{P}(d) = \pi d/2$ and*

$$\mathbb{P}\{w \le \rho\} = \begin{cases} 0, & \rho \le a_3, \\ \dfrac{2\mathcal{P}(\rho)}{\pi d}, & a_3 < \rho < d, \\ 1, & \rho \ge d. \end{cases}$$

Two features of the derivation confine the formula to three terms. The level sets of $w$ on the sphere are circles, and the co-area factor cancels between the radius of the circle and the magnitude of the gradient, leaving an angular measure; and the inclusion-exclusion over the three sign constraints terminates at second order, since a triple overlap would force $w \leq 0$.

The density is not smooth, and its failures of smoothness are classified completely. They occur at the critical values of $w$ on the strata of the arrangement of the great circles $\{u_i = 0\}$, which by Lemma 2.3 are the edge lengths, the face diagonals and the space diagonal, and the type of failure is determined by the dimension of the stratum.

**Theorem B** (Theorem 3.9). *The density is real-analytic away from those values. At an interior edge length it has a corner, with $f'$ jumping by $\frac{2}{\pi a_i a_j}$; at a face diagonal it has a square-root fold of even multiplicity $2 \cdot \#\{j: d_j = d_i\}$, the factor two counting the two level-circle endpoints contributed by each coincident face diagonal; and at the space diagonal a terminal jump of size $4/d$. An edge length and a face diagonal coincide precisely when the edge lengths satisfy $a_1^2 = a_2^2 + a_3^2$, and then only as $a_1 = d_1$, in which case the fold and the corner superpose. The density vanishes below $a_3$, rises linearly from it, increases strictly on $(a_3, \max_i d_i)$, and satisfies $f(\rho) = 4/d$ for $\max_i d_i < \rho < d$.*

The moments are obtained from a generating function rather than term by term: the Gaussian representation of the uniform measure expresses $w$ as a ratio whose numerator is a sum of independent half-normal variables, and the resulting function is entire. A counting rule then governs the powers of $\pi$, the degree of $\kappa_n$ in $\pi^{-1}$ being exactly $\lfloor n/2 \rfloor$, so that $\pi^{-3}$ occurs in no cumulant of order at most five and does occur in the sixth. The mean-normalised cumulants are polynomials in two scale-invariant combinations of the intrinsic volumes, given in closed form through the fifth order.

**Theorem C** (Theorem 6.8, Corollary 6.9). *The first three cumulants determine the edge lengths up to permutation, and hence the entire law:*

$$L = 2\kappa_1, \qquad S = \frac{3\pi}{\pi - 2}\left(\frac{L^2}{12} - \kappa_2\right), \qquad V = \frac{8\pi}{3(4-\pi)}\left(\kappa_3 - \frac{5\pi - 16}{16\pi} LS\right)$$

*after which the edges are the positive roots of $z^3 - Lz^2 + S/2\, z - V$. The initial cumulant sequence is minimal in length: the mean does not determine the variance, and the mean and variance together do not determine the third cumulant. A triple arises from a box precisely when the recovered $L$, $S$ and $V$ satisfy the conditions of Lemma 6.1.*

Finally we ask how far this reaches. The first three moments of the brightness are available for every convex body through the cosine transform of its surface area measure, the third through the classical trivariate Gaussian absolute moment, but the closure in the intrinsic volumes does not extend.

**Theorem D** (Theorem 7.5). *Among convex bodies in $\mathbb{R}^3$, the intrinsic volumes determine neither the variance of the width nor the variance of the brightness. The box with edges $(1,1,4)$ and every right prism of height one over a triangle of area four and perimeter ten share the values $V = 4$, $S = 18$ and $V_1 = 6$; an explicitly constructed member of that prism family has both variances strictly smaller than those of the box.*

The sign-pattern form of the density argument extends to exactly those bodies whose central symmetral $1/2\,(K-K)$ is a zonotope, a class strictly larger than the zonotopes and containing bodies that are not centrally symmetric. For a general parallelepiped we give an exact inner-product test deciding which stratum-critical values are active; we conjecture in Section 7 that the classification of Theorem B extends to all zonotopal central symmetrals. The cellwise co-area argument, by contrast, needs only that $w_K$ be piecewise linear and so applies to every full-dimensional polytope.

Figure 1 shows the density and its derivative for a box whose six edge and face-diagonal values are pairwise distinct.

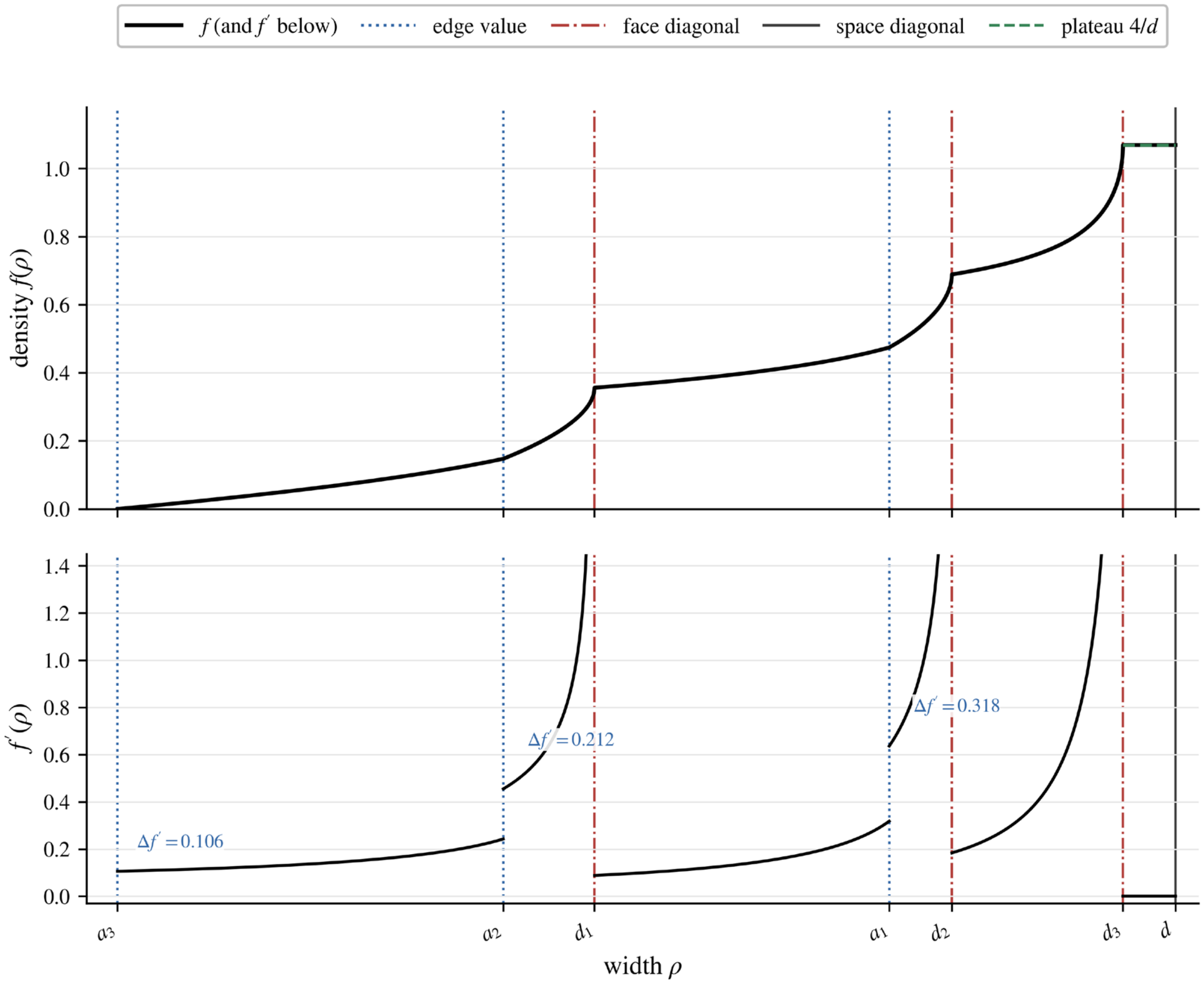

**Figure 1.** The density $f$ of the width of the box with edges $(3,2,1)$, above, and its derivative, below. The support is $[1,\sqrt{14}]$, from the shortest edge to the space diagonal. This box is not edge–face resonant, so its three edge lengths and three face diagonals are distinct and all seven critical values of Lemma 2.3 appear separately. At an edge length the density has a corner, with $f'$ jumping by $\frac{2}{\pi a_i a_j}$ as in Theorem 3.4; the three values are marked. The corner at the shortest edge is the linear vanishing of Corollary 3.5. At a face diagonal the derivative diverges and the density has a square-root fold, by Theorem 3.7. Above the largest face diagonal $f$ is constant at $4/d$, by Proposition 3.3. For a cube the three edge lengths coincide with the lower endpoint of the support and the corners are not separately visible.

## 1.5 Organisation

Section 2 fixes notation and records the elementary geometry of the box: that its elementary symmetric functions are its intrinsic volumes, the two scale-invariant parameters on which the mean-normalised cumulants depend, and the critical values of the width on a natural stratification of the sphere. Section 3 proves Theorems A, A′ and B: the co-area factor for a linear form on the sphere cancels against the density of the level circle, and an inclusion-exclusion over the three excluded arcs terminates at second order, giving the density in closed form; an elementary antiderivative gives the CDF, and the local analysis at each critical value classifies the singularities. Section 4 computes the moments. The Gaussian representation of the uniform measure on $S^2$ expresses $w$ as a ratio whose numerator is a sum of independent half-normal variables; this yields a generating function producing every moment, and with it the counting rule for the powers of $\pi$ on which Section 5 depends. Section 5 proves the degree theorem for the powers of $\pi^{-1}$, establishes the closure of the mean-normalised cumulants in the two scale-invariant parameters, and gives them explicitly through the fifth order. Section 6 proves that the dependence on the intrinsic volumes is graded, that each step of the grading is sharp, and that the first three cumulants determine the box, which is Theorem C. Section 7 examines how far the construction reaches, proving Theorem D, giving an exact incidence criterion for parallelepipeds and stating the conjecture above. Section 8 discusses extensions to zonotopes, general polytopes and higher dimensions, and notes that the cancellation underlying Theorem A is particular to three dimensions. Appendix A records the elementary algebra behind the moment and cumulant formulas.

After the common preliminaries of Section 2, the density analysis of Section 3 and the moment analysis of Sections 4 to 6 are logically independent and may be read in either order.

# 2. Preliminaries

## 2.1 The box and its width

Throughout, $K$ denotes the rectangular box with edge lengths $a_1, a_2, a_3 > 0$, positioned with its edges parallel to the coordinate axes and centred at the origin. Its support function is $h_K(u) = 1/2 \sum_i a_i \, |u_i|$, so that by the definition in Section 1.1

$$w(u) := w_K(u) = \sum_{i=1}^{3} a_i \; |u_i|, \qquad u \in S^2$$

We study the distribution of $w$ when $u$ is uniform on $S^2$. The function $w$ is symmetric in the $a_i$ and homogeneous of degree one in them, so we may assume $a_1 \geq a_2 \geq a_3$ and normalise the scale wherever it is convenient.

Since $w_{cK} = c \, w_K$ pointwise for $c > 0$, the width law scales as

$$f_{cK}(\rho) = \frac{1}{c} \, f_K\left(\frac{\rho}{c}\right), \qquad \kappa_n(cK) = c^n \kappa_n(K)$$

so that the mean-normalised cumulants $g_n = \kappa_n/\kappa_1^n$ are scale invariant. This is the normalisation used throughout.

## 2.2 Intrinsic volumes and elementary symmetric functions

Let $\sigma_1, \sigma_2, \sigma_3$ denote the elementary symmetric polynomials in $a_1, a_2, a_3$, and let $V$, $S$ and $L$ denote the volume of $K$, its surface area, and the sum $a_1 + a_2 + a_3$ of its edge lengths. Then

$$\sigma_1 = L, \qquad \sigma_2 = S/2\,, \qquad \sigma_3 = V$$

The Steiner formula [12, §4.2] for a box reads $|K + \varepsilon B| = V + S\varepsilon + \pi L\varepsilon^2 + 4\pi/3\,\varepsilon^3$, the quadratic term arising from the twelve quarter-cylinders along the edges. Comparison with $\sum_{j=0}^{3} \omega_{3-j}\, V_j(K)\varepsilon^{3-j}$, where $\omega_m$ denotes the volume of the unit $m$-ball, gives

$$V_3(K) = V, \qquad V_2(K) = S/2\,, \qquad V_1(K) = L, \qquad V_0(K) = 1$$

For a box the elementary symmetric functions of the edge lengths therefore coincide with the intrinsic volumes, and every symmetric polynomial in the edges is a polynomial in $V_1, V_2, V_3$. This observation underlies Section 5. In particular the mean width is $\mathbb{E}[w] = \sum_i a_i\ \mathbb{E}|u_i| = L/2 = V_1/2$.

**Lemma 2.1.** *The triple $(L, S, V)$ determines $a_1, a_2, a_3$ up to permutation, namely as the roots of $z^3 - Lz^2 + S/2\,z - V$.*

*Proof.* Immediate from the identification of the $\sigma_k$ above. ▫

**Lemma 2.2 (Newton's inequalities; see [12, §7.2]).** *Let $p_k = \sigma_k/\binom{3}{k}$. Then $p_k^2 \geq p_{k-1}p_{k+1}$ for $k = 1{,}2$, that is*

$$S \leq 2/3\,L^2, \qquad\qquad S^2 \geq 12\,LV$$

*with equality in either exactly when $a_1 = a_2 = a_3$.*

## 2.3 Scale-invariant parameters

Since $w$ is homogeneous of degree one in the edge lengths, its mean-normalised cumulants depend on $K$ only through two scale-invariant parameters. We take these to be

$$G = 1 - \frac{S}{L^2}, \qquad\qquad J = \frac{V}{L^3}$$

Let $d$ denote the length of the space diagonal of $K$. Then $d^2 = \sum_i a_i^2 = L^2 - S$, so that

$$G = (\frac{d}{L})^2$$

By Lemma 2.2 and the positivity of $S$,

$$1/3 \leq G < 1, \qquad \text{equivalently} \qquad 1/\sqrt{3}\ L \leq d < L$$

with equality on the left precisely for the cube. The second inequality of Lemma 2.2 gives $J \leq 1/12\,(1 - G)^2$, again with equality only for the cube.

These two bounds do not delimit the admissible set exactly. Since $(G, J)$ is a function of two shape parameters, its range is a two-dimensional region. Its curved boundary is the image of the locus on which two edge lengths coincide; parametrising that locus as $a = (1,1,t)$ with $t > 0$ gives

$$G(t) = 1 - \frac{2+4t}{(2+t)^2}, \qquad J(t) = \frac{t}{(2+t)^3}$$

with $G'(t) = 4(t-1)/(2+t)^3$ and $J'(t) = 2(1-t)/(2+t)^4$, so that both vanish at $t = 1$. The branches $t < 1$ and $t > 1$ correspond to a repeated longest and a repeated shortest edge respectively. The closure of the region also contains the segment $J = 0$, $1/2 \le G \le 1$: the values $1/2 \le G < 1$ arise as one edge tends to zero, and the endpoint $G = 1$ as two do. This is Remark 6.2, and it is not part of the two-equal-edges locus.

The vanishing of $G'$ and $J'$ at $t = 1$ does not by itself make the cube a cusp, since a parametrisation may be singular where its image is smooth. Writing $t = 1 + x$ and expanding,

$$G = 1/3 + 2/27\, x^2 - 4/81\, x^3 + O(x^4), \qquad J = 1/27 - 1/81\, x^2 + 8/729\, x^3 + O(x^4)$$

The quadratic terms are proportional, so the combination in which they cancel isolates the leading behaviour: with

$$u = G - 1/3\,, \qquad v = J - 1/27 + 1/6\, u$$

one has $u = 2/27\, x^2\big(1 + O(x)\big)$ and $v = 2/729\, x^3\big(1 + O(x)\big)$, whence

$$v^2 = 1/54\ u^3\big(1 + O(x)\big)$$

The two branches $x < 0$ and $x > 0$ are the two halves of this semicubical cusp, meeting at the cube with common tangent direction $u$, that is with slope $-1/6$ in the $(G, J)$ plane.

## 2.4 Critical values of the width

The three great circles $\{u_i = 0\}$ stratify $S^2$ into six 0-strata $\pm e_i$, twelve 1-strata and eight 2-strata. On each stratum the signs of the coordinates are constant, so $w$ restricts to a linear form and is smooth.

**Lemma 2.3.** *The critical values of* $w$ *on the strata of this arrangement are*

$$a_i \quad \text{(0-strata)}, \qquad d_i := \sqrt{a_j^2 + a_k^2} = \sqrt{d^2 - a_i^2} \quad \text{(1-strata)}, \qquad d = \Big(\sum_i a_i^2\Big)^{1/2} \quad \text{(2-strata)}$$

*namely the edge lengths, the face diagonals and the space diagonal.*

*Proof.* On a 2-stratum $w = a \cdot u$ with fixed signs; its only critical point on $S^2$ is $u \propto a$, with value $|a| = d$. A 1-stratum lies in some $\{u_i = 0\}$, on which $w$ restricts to a linear form in the two remaining coordinates, critical where $u$ is proportional to $\big(a_j, a_k\big)$ within that great circle, with value $\sqrt{a_j^2 + a_k^2}$. A 0-stratum is the point $\pm e_i$, where $w = a_i$. ▫

**Remark 2.4.** For an arbitrary convex body, $\max_u w_K(u) = \max_{x,y\in K}|x-y| = \mathrm{diam}(K)$ and $\min_u w_K(u)$ is the thickness of $K$, so that

$$\mathrm{supp}(w_K) = [\,\mathrm{thickness}(K),\ \mathrm{diam}(K)\,]$$

For a box these endpoints are the shortest edge and the space diagonal, recovered as the extreme critical values of Lemma 2.3. Section 3 shows that the density of $w$ is real-analytic away from the values listed there and fails to be real-analytic at each of them, the nature of the failure being determined by the dimension of the stratum on which the value is attained.

# 3. The density and its singularity structure

## 3.1 The density

Throughout this section $d = \left(\sum_i a_i^2\right)^{1/2}$ and $d_i = \left(d^2 - a_i^2\right)^{1/2}$ as in Lemma 2.3, and for $\rho \in (\min_i a_i,\ d)$ we write

$$s = \sqrt{d^2-\rho^2}, \qquad t_i = \frac{\rho\, a_i}{s\, d_i}, \qquad \eta_i = \arccos\bigl(\min(t_i,1)\bigr), \qquad \Theta_{ij} = \arccos\left(\frac{-a_i a_j}{d_i d_j}\right)$$

Since $\rho > 0$ we have $t_i \geq 0$, hence $\eta_i \in [0,\pi/2]$. Throughout, $(x)_+$ denotes $\max(x,0)$.

**Theorem 3.1.** *The distribution of* $w$ *is absolutely continuous, with density*

$$f(\rho) = \frac{2}{\pi d}\,\Psi(\rho), \qquad \Psi(\rho) = 2\pi - 2\sum_{i=1}^{3}\eta_i + \sum_{i<j}(\eta_i+\eta_j-\Theta_{ij})_+ \tag{2}$$

*for* $\rho \in (\min_i a_i, d)$*, and* $f = 0$ *outside* $[\min_i a_i, d]$.

*Proof.* **Reduction to one orthant.** The group of sign changes $u_i \mapsto \pm u_i$ preserves $w$ and the uniform measure on $S^2$, and acts simply transitively on the eight open orthants. Writing

$$O = \{u \in S^2 \colon u_i > 0 \text{ for all } i\}$$

and $\sigma$ for surface measure, the boundaries of the orthants have $\sigma$-measure zero and, for every Borel set $B \subset \mathbb{R}$,

$$\Pr[w \in B] = \frac{8}{4\pi}\,\sigma(\{u \in O \colon w(u) \in B\}) = \frac{2}{\pi}\,\sigma(\{u \in O \colon w(u) \in B\})$$

**Geometry of a level set.** On $O$ the signs are fixed and

$$w(u) = a \cdot u$$

The intrinsic gradient of this restriction to $S^2$ is

$$\nabla_{S^2} w(u) = a - (a\cdot u)u$$

and on the level $w = \rho$ its magnitude is

$$|\nabla_{S^2} w| = \sqrt{|a|^2 - \rho^2} = s$$

The level set is the intersection of $S^2$ with the plane $a \cdot u = \rho$, hence a circle of Euclidean radius $s/d$. If $\Psi(\rho)$ denotes the angular measure of the part of that circle lying in $O$, its one-dimensional Hausdorff measure is

$$\mathcal{H}^1(\{w = \rho\} \cap O) = \frac{s}{d}\,\Psi(\rho)$$

Consequently

$$\int_{\{w=\rho\}\cap O} \frac{d\mathcal{H}^1}{|\nabla_{S^2} w|} = \frac{1}{s} \cdot \frac{s}{d}\,\Psi(\rho) = \frac{\Psi(\rho)}{d}$$

The radius of the level circle and the magnitude of the intrinsic gradient cancel.

**Computation of $\Psi$.** Let $n = a/d$ and let $P$ denote orthogonal projection onto $n^{\perp}$. Parametrise the level circle as

$$u = \frac{\rho}{d} n + \frac{s}{d} v$$

with $v$ a unit vector of $n^{\perp}$. Since

$$|Pe_i| = \sqrt{1 - a_i^2/d^2} = d_i/d$$

write $\hat{e}_i = Pe_i/|Pe_i|$, and let $v$ and $\hat{e}_i$ have angular coordinates $\varphi$ and $\varphi_i$ in $n^{\perp}$. Then

$$u_i = \frac{1}{d^2}[\rho a_i + s d_i \cos(\varphi - \varphi_i)]$$

Thus $u_i \leq 0$ precisely when

$$\cos(\varphi - \varphi_i) \leq -t_i$$

where $t_i = \rho a_i/(s d_i)$. With $\eta_i = \arccos\big(\min(t_i, 1)\big) \in [0, \pi/2]$, this is the arc $A_i$ centred at $\varphi_i + \pi$ with half-width $\eta_i$ and length $2\eta_i$; the arc is empty when $t_i \geq 1$.

The separation between the centres of $A_i$ and $A_j$ is the angle $\Theta_{ij}$, because

$$\frac{Pe_i \cdot Pe_j}{|Pe_i||Pe_j|} = \frac{-a_i a_j}{d_i d_j} = \cos\Theta_{ij}$$

Since $\eta_i, \eta_j \leq \pi/2$,

$$\eta_i + \eta_j \leq \pi \leq 2\pi - \Theta_{ij}$$

and therefore $A_i$ and $A_j$ meet in at most one arc, of length

$$|A_i \cap A_j| = (\eta_i + \eta_j - \Theta_{ij})_+$$

The triple intersection is empty. Indeed, a point of $A_1 \cap A_2 \cap A_3$ would satisfy $u_i \leq 0$ for every $i$ and hence

$$\rho = w(u) = \sum_i a_i u_i \leq 0$$

contrary to $\rho > 0$. Inclusion-exclusion therefore terminates at second order and gives

$$\Psi(\rho) = 2\pi - |A_1 \cup A_2 \cup A_3| = 2\pi - 2\sum_i \eta_i + \sum_{i<j} (\eta_i + \eta_j - \Theta_{ij})_+$$

**Absolute continuity.** Let $B$ be a Borel subset of $(\min_i a_i, d)$. On $w^{-1}(B) \cap O$ the intrinsic gradient does not vanish. The integral form of the co-area formula, applied to the indicator of this set divided by $|\nabla_{S^2} w|$, gives

$$\begin{aligned} \sigma(\{u \in O: w(u) \in B\}) &= \int_B \left( \int_{\{w=\rho\}\cap O} \frac{d\mathcal{H}^1}{|\nabla_{S^2} w|} \right) d\rho \\ &= \int_B \frac{\Psi(\rho)}{d}\, d\rho \end{aligned}$$

Combining this identity with the reduction to one orthant yields

$$\Pr[w \in B] = \int_B \frac{2\Psi(\rho)}{\pi d}\, d\rho$$

for every such Borel set $B$. The endpoint level sets are finite unions of arcs or points and have spherical area zero, so they carry no mass. By Remark 2.4 the range of $w$ is $[\min_i a_i, d]$. Hence the law is absolutely continuous and has the stated density, with the displayed representative taken to be zero outside that interval. In particular, the terminal discontinuity of this representative at $d$ is not an atom of the distribution. ▫

**Proposition 3.2.** $\Theta_{12} + \Theta_{23} + \Theta_{31} = 2\pi$. Consequently, if all three positive parts in (2) are strictly positive then $\Psi(\rho) = 0$.

*Proof.* $\sum_i a_i\, Pe_i = P(a) = P(dn) = 0$, so the three nonzero coplanar vectors $Pe_i$ admit a vanishing positive combination and therefore do not lie in a closed half-plane. Their directions divide the circle into three consecutive gaps summing to $2\pi$, each less than $\pi$; for three directions each pairwise angle in $[0, \pi]$ is one of those gaps. If every positive part is active then

$$\Psi = 2\pi - 2\sum_i \eta_i + \sum_{i<j} (\eta_i + \eta_j - \Theta_{ij}) = 2\pi - 2\sum_i \eta_i + 2\sum_i \eta_i - 2\pi = 0$$

▫

## 3.2 Elementary structure

**Proposition 3.3.** Let $a_1 \geq a_2 \geq a_3$. Then

1. $f \equiv 4/d$ on $(\max_i d_i,\ d) = \left(\sqrt{a_1^2 + a_2^2},\ d\right)$;
2. $a_3 < \min_i d_i$, and $a_k < \max_i d_i$ for every $k$;
3. an edge length equals a face diagonal if and only if $a_1^2 = a_2^2 + a_3^2$, in which case $a_1 = d_1$, and this is the only equality between an edge length and a face diagonal.

*Proof.* (1) For $\rho > d_i$ one has $t_i > 1$ and $\eta_i = 0$; if all three vanish then every positive part vanishes and $\Psi = 2\pi$. (2) $\min_i d_i = \sqrt{a_2^2 + a_3^2} > a_3$ and $\max_i d_i = \sqrt{a_1^2 + a_2^2} > a_1 \geq a_k$. (3) If $a_k = d_i$ with $k \neq i$ then $a_k^2 + a_i^2 = d^2$, forcing the remaining edge to vanish. If $a_i = d_i$ then $a_i^2 = a_j^2 + a_k^2$, so $a_i$ exceeds both $a_j$ and $a_k$ and therefore $i = 1$. ▫

A box satisfying $a_1^2 = a_2^2 + a_3^2$, equivalently $a_1 = d_1$, will be called **edge–face resonant**; examples are $(\sqrt{2}, 1, 1)$, (5,4,3) and (13,12,5).

By (2) the least critical value is $a_3$ and the greatest is $d$, and the remaining values interleave in an order depending on the box. By (3) an edge length and a face diagonal coincide only in the resonant case, and then only as $a_1 = d_1$.

### 3.3 Corners at the edge lengths

**Theorem 3.4.** Let $\{i, j, k\} = \{1,2,3\}$.

1. The term $\left(\eta_i + \eta_j - \Theta_{ij}\right)_+$ is positive for $\rho < a_k$ and vanishes for $\rho \geq a_k$.
2. If moreover $a_k \neq d_k$, so that by Proposition 3.3(3) the value $a_k$ is not a face diagonal, then $f$ has finite one-sided derivatives at $a_k$ and

$$f'(a_k^+) - f'(a_k^-) = \frac{2}{\pi\, a_i a_j} \tag{3}$$

If several edge lengths coincide the corresponding jumps add.

*Proof.* Put $\rho = a_k$, so that $s = d_k$ and $\cos\eta_i = a_k a_i/(d_k d_i)$. From

$$d_k^2 d_i^2 - a_k^2 a_i^2 = \left(d^2 - a_k^2\right)\left(d^2 - a_i^2\right) - a_k^2 a_i^2 = d^2\left(d^2 - a_i^2 - a_k^2\right) = d^2 a_j^2 \tag{4}$$

we obtain $\sin\eta_i = d\, a_j/(d_k d_i)$, and by (4) with $i$ and $j$ interchanged, $\sin\eta_j = d\, a_i/\left(d_k d_j\right)$. Hence

$$\cos\left(\eta_i + \eta_j\right) = \frac{a_k a_i}{d_k d_i} \cdot \frac{a_k a_j}{d_k d_j} - \frac{d\, a_j}{d_k d_i} \cdot \frac{d\, a_i}{d_k d_j} = \frac{a_i a_j\left(a_k^2 - d^2\right)}{d_k^2 d_i d_j} = \frac{-a_i a_j}{d_i d_j} = \cos\Theta_{ij}$$

Both $\eta_i + \eta_j$ and $\Theta_{ij}$ lie in $[0, \pi]$, on which the cosine is injective, so $\eta_i + \eta_j = \Theta_{ij}$ at $\rho = a_k$.

Differentiating $\eta_i = \arccos t_i$ gives $\eta_i{}' = -t_i{}'/\sqrt{1 - t_i^2}$, and from $t_i = (a_i/d_i)\rho\,(d^2 - \rho^2)^{-1/2}$ we get $t_i{}' = (a_i/d_i)\, d^2 (d^2 - \rho^2)^{-3/2}$. At $\rho = a_k$,

$$t_i{}' = \frac{a_i\, d^2}{d_i\, d_k^3}, \qquad \sqrt{1 - t_i^2} = \frac{d\, a_j}{d_k d_i}, \qquad \text{so} \qquad \eta_i{}'(a_k) = -\frac{a_i\, d}{a_j\, d_k^2}$$

Adding the same expression with $i$ and $j$ interchanged and using $a_i^2 + a_j^2 = d_k^2$,

$$\eta_i{}'(a_k) + \eta_j{}'(a_k) = -\frac{d}{d_k^2} \cdot \frac{a_i^2 + a_j^2}{a_i a_j} = -\frac{d}{a_i a_j}$$

The argument $\eta_i + \eta_j - \Theta_{ij}$ therefore vanishes at $a_k$ with strictly negative derivative. To pass from this local statement to the global one, note that $t_\ell{}' = (a_\ell/d_\ell)\, d^2(d^2 - \rho^2)^{-3/2} > 0$ wherever $t_\ell < 1$. Thus every $\eta_\ell$ is non-increasing, and is strictly decreasing while $\rho < d_\ell$. Moreover $a_k < d_i$ and $a_k < d_j$, so $\eta_i + \eta_j - \Theta_{ij}$ is strictly decreasing throughout $\rho < a_k$. It therefore crosses zero exactly once, at $a_k$, and cannot become positive again, which gives assertion 1. For assertion 2, the hypothesis $a_k \neq d_k$ ensures that every $\eta_\ell$ is differentiable at $a_k$, so the term contributes $\eta_i{}' + \eta_j{}'$ to $\Psi'$ below $a_k$ and nothing above, whence $\Psi'(a_k^+) - \Psi'(a_k^-) = d/(a_i a_j)$; multiplying by $\frac{2}{\pi d}$ gives the jump. Additivity under coincidence is immediate, distinct pairs contributing distinct terms. ▫

**Corollary 3.5.** At the lower endpoint $\rho = a_3 = \min_i a_i$ the density vanishes and rises linearly, with

$$f'(a_3^+) = \sum_{k\,:\,a_k = a_3} \frac{2}{\pi\, a_i a_j}$$

In particular this slope equals $\frac{2}{\pi a_1 a_2}$ when the minimum is simple, $\frac{4}{\pi a_1 a_3}$ when it has multiplicity two, and $\frac{6}{\pi a^2}$ for the cube of edge $a$.

*Proof.* $f = 0$ below $a_3$, so the left derivative vanishes and Theorem 3.4 supplies the jump. ▫

**Theorem 3.6 (monotonicity).** The density $f$ is strictly increasing on $(\min_i a_i,\ \max_i d_i)$ and constant, equal to $4/d$, on $(\max_i d_i,\ d)$.

*Proof.* On an interval free of critical values, let $\mathcal{J}$ be the set of pairs $\{i,j\}$ whose term in $\Psi$ is positive, and put $r_i = \#\{\{i,j\} \in \mathcal{J}\}$, so $r_i \in \{0,1,2\}$. Differentiating

$$\Psi = 2\pi - 2\sum_i \eta_i + \sum_{\{i,j\} \in \mathcal{J}} \left(\eta_i + \eta_j - \Theta_{ij}\right)$$

collects the coefficient of $\eta_i{}'$ as $-2 + r_i$, so that

$$\Psi'(\rho) = \sum_{i=1}^{3} (r_i(\rho) - 2)\ \eta_i{}'(\rho)$$

Where $t_i < 1$ one has $t_i{}' = (a_i/d_i) d^2 (d^2 - \rho^2)^{-3/2} > 0$ and hence $\eta_i{}' = -t_i{}'/\sqrt{1 - t_i^2} < 0$; where $t_i \geq 1$ one has $\eta_i \equiv 0$ and $\eta_i{}' = 0$. Since $r_i - 2 \leq 0$ in every case, each term is non-negative and $f' \geq 0$.

For strictness, note that by Theorem 3.4(1) the pair $\{i,j\}$ is active exactly when $\rho < a_k$, so $r_i = 2$ requires $\rho < \min(a_j, a_k)$. Order $a_1 \geq a_2 \geq a_3$ and take $\rho > a_3$. Then $r_1 \leq 1$ and $r_2 \leq 1$, while $\eta_1 > 0$ for $\rho < d_1$ and $\eta_2 > 0$ for $\rho < d_2$; and for $\rho > a_2$ also $r_3 \leq 1$, with $\eta_3 > 0$ for $\rho < d_3$. Since $a_3 < d_1 \leq d_2 \leq d_3$ and $a_2 < d_2$, some index contributes a strictly positive term at every $\rho$ in $(a_3, d_3)$, whence $f' > 0$ there. Above $\max_i d_i$ every $t_i$ exceeds one, so $\Psi = 2\pi$ and $f = 4/d$ by Proposition 3.3(1). ▫

Since $f > 0$ on $(\min_i a_i, d)$, at most two of the three positive parts can be active at any point of the support: if all three were active, Proposition 3.2 would give $\Psi = 0$ and hence $f = 0$. Thus simultaneous activity of all three can occur only off the support, where the displayed density representative vanishes.

## 3.4 Folds at the face diagonals

**Theorem 3.7.** Let $\beta_i = d^2/(d_i\, a_i^2)$. As $\rho \uparrow d_i$,

$$f(d_i) - f(\rho) = \frac{2\nu_i}{\pi d}\sqrt{2\beta_i}\,\sqrt{d_i - \rho}\,\big(1 + o(1)\big), \qquad \nu_i = 2 \cdot \#\{j : d_j = d_i\} \tag{5}$$

---

In particular $\nu_i$ is even: each face diagonal $d_j$ equal to $d_i$ contributes the two endpoints of its excluded arc, so $\nu_i$ counts endpoints rather than diagonals and a simple face diagonal has $\nu_i = 2$.

*Proof.* At $\rho = d_i$ one has $s = a_i$ and $t_i = 1$. From $t_i' = (a_i/d_i)\, d^2(d^2 - \rho^2)^{-3/2}$, evaluation at $\rho = d_i$ gives $t_i'(d_i) = d^2/(d_i a_i^2) = \beta_i$, so $t_i = 1 - \beta_i(d_i - \rho) + O((d_i - \rho)^2)$. Since $\arccos(1 - \varepsilon) = \sqrt{2\varepsilon}\,\big(1 + O(\varepsilon)\big)$,

$$\eta_i = \sqrt{2\beta_i(d_i - \rho)}\,\big(1 + o(1)\big)$$

The only terms of $\Psi$ that are not real-analytic at $d_i$ are those containing $\eta_i$, namely $-2\eta_i$ and the two positive parts indexed by pairs containing $i$. By Theorem 3.4(1) the pair $\{i,j\}$ contributes only for $\rho < a_k$, and $a_k < \sqrt{a_j^2 + a_k^2} = d_i$; likewise the pair $\{i,k\}$ contributes only for $\rho < a_j < d_i$. Both are therefore inactive in a neighbourhood of $d_i$, and $\eta_i$ enters $\Psi$ only through $-2\eta_i$, contributing $2\eta_i$ to $\Psi(d_i) - \Psi(\rho)$.

If the box is edge–face resonant and $i = 1$, so that $d_1 = a_1$, the remaining pair $\{2,3\}$ is not inactive: by Theorem 3.4 its argument vanishes at $a_1$, and by the derivative computed there

$$\eta_2 + \eta_3 - \Theta_{23} = \frac{d}{a_2 a_3}\,(a_1 - \rho) + O((a_1 - \rho)^2) \qquad (\rho \uparrow a_1)$$

so that term is $O(a_1 - \rho)$ and is subordinate to $\eta_1$, which is of order $\sqrt{d_1 - \rho}$. The leading expansion is therefore unaffected.

If $\#\{j : d_j = d_i\} = m$ then $m$ such terms are simultaneously singular, and $d_j = d_i$ forces $a_j = a_i$ and hence $\beta_j = \beta_i$. Multiplying by $\frac{2}{\pi d}$ gives the stated form with $\nu_i = 2m$, which is even and at least 2. ▫

## 3.5 Superposition at a resonance

**Theorem 3.8.** Let the box be edge–face resonant, so that $a_1 = d_1$, and put $C_1 = \frac{2\nu_1}{\pi d}\sqrt{2\beta_1}$ in the notation of Theorem 3.7. Then

$$f(a_1) - f(\rho) = C_1\sqrt{a_1 - \rho} + O(a_1 - \rho) \qquad (\rho \uparrow a_1)$$

so that $f'(\rho) \to +\infty$ as $\rho \uparrow a_1$ while $f'(a_1^+)$ is finite; and the desingularised density

$$\tilde{f}(\rho) = f(\rho) + C_1\sqrt{(a_1 - \rho)_+}$$

has finite one-sided derivatives at $a_1$, with

$$\tilde{f}\,'(a_1^+) - \tilde{f}\,'(a_1^-) = \frac{2}{\pi\, a_2 a_3}$$

*Proof.* Put $\delta = a_1 - \rho$. Since $a_1 = d_1$ we have $t_1 = 1$ at $\rho = a_1$, while $t_2 = a_2/d_2 < 1$ and $t_3 = a_3/d_3 < 1$ there, so $\eta_2$ and $\eta_3$ are real-analytic in a neighbourhood of $a_1$. By Theorem 3.4(1) the terms $(\eta_1 + \eta_2 - \Theta_{12})_+$ and $(\eta_1 + \eta_3 - \Theta_{13})_+$ vanish for $\rho \geq a_3$ and $\rho \geq a_2$ respectively, hence identically near $a_1$; and $(\eta_2 + \eta_3 - \Theta_{23})_+$ vanishes for $\rho \geq a_1$. Writing

$$R(\rho) = \frac{2}{\pi d}\big(2\pi - 2\eta_2(\rho) - 2\eta_3(\rho)\big)$$

which is real-analytic at $a_1$, we therefore have

$$f(\rho) = R(\rho) - \frac{4}{\pi d}\,\eta_1(\rho) + \frac{2}{\pi d}(\eta_2 + \eta_3 - \Theta_{23})_+$$

Since $d^2 - a_1^2 = d_1^2 = a_1^2$, the derivative $t_1{}' = (a_1/d_1)d^2(d^2 - \rho^2)^{-3/2}$ takes the value $d^2/(d_1 a_1^2) = \beta_1$ at $\rho = a_1$; as $t_1$ is real-analytic there with $t_1(a_1) = 1$, we have $1 - t_1 = \beta_1\delta\big(1 + O(\delta)\big)$ with analytic remainder. Combining this with $\arccos(1 - \varepsilon) = \sqrt{2\varepsilon}\,h(\varepsilon)$, where $h$ is analytic near 0 with $h(0) = 1$, gives

$$\eta_1 = \delta^{1/2}A(\delta) \;\; (\delta \geq 0), \qquad \eta_1 \equiv 0 \;\; (\delta \leq 0)$$

with $A$ analytic near 0 and $A(0) = \sqrt{2\beta_1}$. By Theorem 3.4 the argument $\eta_2 + \eta_3 - \Theta_{23}$ is analytic and vanishes at $a_1$ with $\eta_2{}'(a_1) + \eta_3{}'(a_1) = -d/(a_2 a_3)$, so

$$(\eta_2 + \eta_3 - \Theta_{23})_+ = \frac{d}{a_2 a_3}\,\delta + O(\delta^2) \;\; (\delta \geq 0), \qquad \equiv 0 \;\; (\delta \leq 0)$$

Resonance gives $d_1 < d_2 \leq d_3$ and hence $\nu_1 = 2$, so $C_1 = 4/\pi d\sqrt{2\beta_1} = 4/\pi d\, A(0)$. For $\delta \geq 0$,

$$\tilde{f}(\rho) = R(\rho) + \frac{4}{\pi d}\,\delta^{1/2}\big(A(0) - A(\delta)\big) + \frac{2}{\pi d}(\eta_2 + \eta_3 - \Theta_{23})_+$$

The middle term is $\delta^{3/2}$ times an analytic function, hence differentiable at $\delta = 0$ with vanishing derivative. Therefore $\tilde{f}\,'(a_1^-) = R'(a_1) - 2/\pi a_2 a_3$. For $\delta \le 0$ all three added terms vanish and $\tilde{f} = R$, so $\tilde{f}\,'(a_1^+) = R'(a_1)$. Subtracting gives the stated jump. ▫

Thus at a resonance the two mechanisms superpose: the square-root fold dominates the visible singularity, and the edge corner survives in the regular remainder.

## 3.6 Classification

**Theorem 3.9 (complete singularity structure).** Let $a_1 \ge a_2 \ge a_3 > 0$. The density $f$, in the representative furnished by Theorem 3.1, is real-analytic on every component of

$$a_3 < \rho < d, \qquad \rho \ne a_1,\ \rho \ne a_2,\ \rho \ne d_1,\ \rho \ne d_2,\ \rho \ne d_3$$

coincident critical values being counted once. Its behaviour at the exceptional values is as follows.

1. **Corner.** At an interior edge value $\rho = a_k > a_3$ with $a_k \ne d_k$ the one-sided derivatives are finite and

$$f'(a_k^+) - f'(a_k^-) = \frac{2}{\pi\, a_i a_j}, \qquad \{i,j,k\} = \{1,2,3\}$$

contributions adding when several interior edge lengths coincide.

2. **Fold.** At $\rho = d_i$ with $d_i \ne a_i$, writing $\beta_i = d^2/\left(d_i a_i^2\right)$,

$$f(d_i) - f(\rho) = \frac{2\nu_i}{\pi d}\sqrt{2\beta_i}\,\sqrt{d_i - \rho}\,\left(1 + o(1)\right), \qquad \nu_i = 2 \cdot \#\{j\colon d_j = d_i\}$$

where the factor two in $\nu_i$ counts the two level-circle endpoints contributed by each coincident face diagonal. Thus $\nu_i$ is even; moreover, $f$ is right-analytic at $d_i$ and $f'(\rho) \to +\infty$ as $\rho \uparrow d_i$.

3. **Superposition.** The only possible edge–face coincidence is $a_1 = d_1$. There both mechanisms act: $f'(\rho) \to +\infty$ as $\rho \uparrow a_1$, whereas $f'(a_1^+)$ is finite. The resonant face diagonal is simple, since $d_1 = d_2$ would force $a_1 = a_2$ and $d_1 = d_3$ would force $a_1 = a_3$, either of which makes a third edge vanish; hence $\nu_1 = 2$. Putting

$$C_{\mathrm{f}} := \frac{4}{\pi d}\sqrt{2\beta_1}, \qquad \tilde{f}(\rho) := f(\rho) + C_{\mathrm{f}}\sqrt{(a_1 - \rho)_+}$$

the desingularised density retains the corner:

$$\tilde{f}\,'(a_1^+) - \tilde{f}\,'(a_1^-) = \frac{2}{\pi\, a_2 a_3}$$

4. **Lower endpoint.** The density vanishes on $(-\infty, a_3]$ and rises linearly from the endpoint, with

$$f'(a_3^+) = \sum_{k\,:\,a_k = a_3} \frac{2}{\pi\, a_i a_j}, \qquad \{i,j,k\} = \{1,2,3\}$$

5. **Plateau and terminal step.** One has $f \equiv 4/d$ on $(\max_i d_i,\ d)$ and $f = 0$ above $d$, so that

$$\lim_{\rho\uparrow d} f(\rho) = \frac{4}{d}, \qquad \lim_{\rho\downarrow d} f(\rho) = 0$$

a terminal jump of size $4/d$.

Moreover $f$ is strictly increasing on $(a_3,\ \max_i d_i)$.

*Proof.* Away from the listed values the set of active positive-part terms is locally constant, and each $\eta_i$ is either identically zero or equals $\arccos t_i$ with $t_i \in (0,1)$ real-analytic in $\rho$; since arccos is real-analytic on $(-1,1)$, so is $f$. Assertion 1 is Theorem 3.4(2), assertion 2 is Theorem 3.7, assertion 3 is Theorem 3.8 with Proposition 3.3(3), assertion 4 is Corollary 3.5, and assertion 5 follows from Proposition 3.3(1) and Theorem 3.1. Strict monotonicity is Theorem 3.6. ▫

**Remark 3.10.** The dimension of the stratum in Lemma 2.3 determines the mechanism: a 0-stratum produces a corner, a 1-stratum a square-root fold, and the 2-strata the terminal step. When a critical value is attained on strata of two dimensions the mechanisms superpose, as in Theorem 3.9(3). For a cube all three edge lengths coincide with the lower endpoint, so their corner contributions are absorbed into Theorem 3.9(4).

### 3.7 The cumulative distribution function

The global density has an elementary antiderivative. Extend $\eta_i$ continuously to $d$ by $\eta_i(d) = 0$, and for $\rho \in [a_3, d]$ define

$$E_i(\rho) = \rho\eta_i(\rho) - d\arctan\left(\frac{\sqrt{\left(d_i^2 - \rho^2\right)_+}}{a_i}\right).$$

For $\{i, j, k\} = \{1,2,3\}$ put

$$E_{ij}(\rho) = \begin{cases} E_i(\rho) + E_j(\rho) - \Theta_{ij}\rho, & \rho \le a_k, \\ E_i(a_k) + E_j(a_k) - \Theta_{ij}a_k, & \rho \ge a_k, \end{cases}$$

and

$$\mathcal{P}(\rho) = 2\pi\rho - 2\sum_{i=1}^{3} E_i\,(\rho) + \sum_{i<j} E_{ij}\,(\rho).$$

**Proposition 3.11 (closed CDF).** *The cumulative distribution function of $w$ is*

$$F(\rho) = \begin{cases} 0, & \rho \le a_3, \\ \dfrac{2\mathcal{P}(\rho)}{\pi d}, & a_3 < \rho < d, \\ 1, & \rho \ge d. \end{cases}$$

*Moreover $\mathcal{P}$ is continuously differentiable on $(a_3, d)$, $\mathcal{P}' = \Psi$ there, and*

$$\mathcal{P}(a_3) = 0, \qquad \mathcal{P}(d) = \frac{\pi d}{2}.$$

*Proof.* For $\rho < d_i$, direct differentiation gives

$$\eta_i{}'(\rho) = -\frac{a_i d}{(d^2 - \rho^2)\sqrt{d_i^2 - \rho^2}}.$$

Differentiating $E_i$ shows that the term $\rho\eta_i{}'$ cancels the derivative of the arctangent, and hence $E_i{}' = \eta_i$. For $\rho \geq d_i$ both sides vanish. By Theorem 3.4(1),

$$\left(\eta_i + \eta_j - \Theta_{ij}\right)_+ = \begin{cases} \eta_i + \eta_j - \Theta_{ij}, & \rho < a_k, \\ 0, & \rho \geq a_k, \end{cases}$$

so $E_{ij}{}' = \left(\eta_i + \eta_j - \Theta_{ij}\right)_+$ and therefore $\mathcal{P}' = \Psi$ away from the critical values. In fact $E_i{}'$ tends to the common value 0 at $d_i$, while at $a_k$ the left derivative of $E_{ij}$ is $\eta_i(a_k) + \eta_j(a_k) - \Theta_{ij} = 0$, equal to its right derivative. Thus $\mathcal{P} \in C^1(a_3, d)$ and $\mathcal{P}' = \Psi$ throughout that interval.

It remains to evaluate the endpoints, which also verifies the normalisation without invoking the total mass. At $\rho = a_3$ every $E_{ij}$ is in its first branch. Each $E_i(a_3)$ then occurs twice in $\sum_{i<j} E_{ij}(a_3)$ and cancels the term $-2\sum_i E_i(a_3)$, while Proposition 3.2 gives $\sum_{i<j} \Theta_{ij} = 2\pi$. Hence $\mathcal{P}(a_3) = 0$.

At $\rho = d$, $E_i(d) = 0$. For $\{i,j,k\} = \{1,2,3\}$, the switching-point identity $\eta_i(a_k) + \eta_j(a_k) = \Theta_{ij}$ from Theorem 3.4 and $d_i^2 - a_k^2 = a_j^2$, $d_j^2 - a_k^2 = a_i^2$ give

$$E_{ij}(d) = -d\left[\arctan\left(\frac{a_j}{a_i}\right) + \arctan\left(\frac{a_i}{a_j}\right)\right] = -\frac{\pi d}{2}.$$

There are three pairs, and therefore $\mathcal{P}(d) = 2\pi d - 3\pi d/2 = \pi d/2$. Integrating $f = 2\mathcal{P}'/(\pi d)$ now proves all three branches of the formula. ▫

Thus exact interval and tail probabilities require only algebraic functions, arccosines and arctangents; quantiles reduce to a one-dimensional inversion of this explicit monotone function.

## 3.8 The planar analogue

**Proposition 3.12.** Let $R$ be a rectangle with sides $a_1 \geq a_2$, let $u$ be uniform on $S^1$, and put $d = \sqrt{a_1^2 + a_2^2}$. Then $w = a_1|u_1| + a_2|u_2|$ has density

$$f^{(2)}(\rho) = \frac{2}{\pi} \cdot \frac{N(\rho)}{\sqrt{d^2 - \rho^2}}, \qquad N(\rho) = \begin{cases} 1, & a_2 < \rho < a_1 \\ 2, & a_1 < \rho < d \end{cases}$$

and consequently

$$\mathbb{E}[w] = \frac{2(a_1 + a_2)}{\pi} = \frac{\mathrm{per}(R)}{\pi}, \qquad \mathbb{E}[w^2] = \frac{a_1^2 + a_2^2}{2} + \frac{2a_1a_2}{\pi}$$

*Proof.* Writing $u = (\cos\varphi, \sin\varphi)$, the map $\varphi \mapsto w$ is invariant under $\varphi \mapsto \varphi + \pi$ and $\varphi \mapsto \pi - \varphi$, so $[0, \pi/2]$ is a fundamental domain carrying effective density $2/\pi$. On it $w = d\cos(\varphi - \varphi_0)$ with $\varphi_0 = \arctan(a_2/a_1)$, increasing from $a_1$ to $d$ on $[0, \varphi_0]$ and decreasing from $d$ to $a_2$ on $[\varphi_0, \pi/2]$. Hence $\rho \in (a_2, a_1)$ has one preimage and $\rho \in (a_1, d)$ has two. Since $|dw/d\varphi| = d|\sin(\varphi - \varphi_0)| = \sqrt{d^2 - w^2}$, the change of variables gives the stated density. The moments follow by the substitution $\rho = d\sin\theta$. ▫

**Remark 3.13.** For the unit square Proposition 3.12 gives $\mathbb{E}[w^2] = 1 + 2/\pi$, the planar value obtained in [4]. Total mass, mean and second moment were checked numerically for the rectangles (1,1), (2,1) and (3,2), with agreement to $10^{-9}$.

# 4. Moments

## 4.1 The Gaussian representation

The same Gaussian polar device appears in the general width-moment identity of Kabluchko, Litvak and Zaporozhets [7, Theorem 3.1], which applies to every compact convex body. The lemma below records its direct three-dimensional box form. We include the radial proof because in the present setting it is elementary and fixes the normalisation needed later.

**Lemma 4.1.** Let $g = (g_1, g_2, g_3)$ have independent standard normal coordinates. Then $u = g/|g|$ is uniform on $S^2$ and is independent of $|g|$, and

$$w = \frac{W}{|g|}, \qquad W := \sum_{i=1}^{3} a_i\,|g_i|$$

Consequently, for every $N \geq 0$,

$$m_N = \frac{\mathbb{E}[W^N]}{\mathbb{E}[|g|^N]}, \qquad \mathbb{E}[|g|^N] = \frac{2^{N/2+1}\,\Gamma(3 + N/2)}{\sqrt{\pi}}$$

*Proof.* The law of $g$ is invariant under $O(3)$, so the law of $u$ is the unique rotation-invariant probability measure on $S^2$, and the pair $(|g|, u)$ is independent because the polar decomposition of the density factorises. Since $|g| > 0$ almost surely, $w(u) = \sum_i a_i\,|u_i| = \sum_i a_i\,|g_i|/|g| = W/|g|$. Taking $N$-th powers and using independence gives $\mathbb{E}[W^N] = \mathbb{E}[|g|^N]\,\mathbb{E}[w^N]$, and $\mathbb{E}[|g|^N]$ is finite and positive because $|g|$ has the chi distribution with three degrees of freedom, whose moments are $2^{N/2}\Gamma(3 + N/2)/\Gamma(3/2)$ with $\Gamma(3/2) = \sqrt{\pi}/2$. ▫

**Corollary 4.2.** For non-negative integers $q_1, q_2, q_3$ with $Q = \sum_i q_i$,

$$\mathbb{E}\prod_{i=1}^{3}|u_i|^{q_i} = \frac{1}{2\pi\,\Gamma(3 + Q/2)}\prod_{i=1}^{3}\Gamma\left(\frac{q_i + 1}{2}\right)$$

*Proof.* By Lemma 4.1 the left side equals $\mathbb{E}\prod_i |g_i|^{q_i}/\mathbb{E}[|g|^Q]$, and the numerator factorises as $\prod_i \mathbb{E}\,|g_i|^{q_i}$ with $\mathbb{E}|g_1|^q = 2^{q/2}\Gamma(q + 1/2)/\sqrt{\pi}$. Dividing by the expression for $\mathbb{E}[|g|^Q]$ gives the result, the powers of 2 cancelling. ▫

## 4.2 A generating function for all moments

**Proposition 4.3.** With $d^2 = \sum_i a_i^2$ as in Lemma 2.3, the random variable $W$ has moment generating function

$$\mathbb{E}[e^{tW}] = \exp(1/2\, t^2 d^2) \prod_{i=1}^{3} \left[1 + \operatorname{erf}\left(t a_i/\sqrt{2}\right)\right], \qquad t \in \mathbb{R}$$

which extends to an entire function of $t$. In particular $W$ has finite moments of every order, and by Lemma 4.1 so does $w$, with

$$m_N = \frac{\sqrt{\pi}}{2^{N/2+1}\,\Gamma(3 + N/2)} \cdot \frac{d^N}{dt^N}\Big|_{t=0} \exp(1/2\, t^2 d^2) \prod_{i=1}^{3} \left[1 + \operatorname{erf}\left(t a_i/\sqrt{2}\right)\right]$$

*Proof.* For a standard normal $Z$ and real $s$, completing the square gives

$$\mathbb{E}\left[e^{s|Z|}\right] = 2\int_0^\infty e^{sy} \frac{e^{-y^2/2}}{\sqrt{2\pi}}\, dy = 2e^{s^2/2} \int_0^\infty \frac{e^{-(y-s)^2/2}}{\sqrt{2\pi}}\, dy = 2e^{s^2/2}\Phi(s) = e^{s^2/2}\left[1 + \operatorname{erf}\left(s/\sqrt{2}\right)\right]$$

using $2\Phi(s) = 1 + \operatorname{erf}\left(s/\sqrt{2}\right)$. The coordinates $g_i$ are independent, so $\mathbb{E}[e^{tW}] = \prod_i \mathbb{E}\left[e^{t a_i |g_i|}\right]$, and collecting the exponentials gives $\exp\left(1/2\, t^2 \sum_i a_i^2\right)$. Both exp and erf are entire, so the product is entire and its Taylor coefficients at the origin are the moments of $W$ divided by factorials. The displayed formula for $m_N$ is then Lemma 4.1. ▫

**Remark 4.4.** Proposition 4.3 gives every moment, not only those tabulated below. For the unit cube it yields, for instance,

$$m_6 = \frac{1012 + 135\pi}{35\pi} = 13.060845852, \qquad m_7 = \frac{129}{8} + \frac{14}{\pi} = 20.581338407$$

The generating function is that of $W$ and not of $w$; the normalisation by $\mathbb{E}[|g|^N]$ in Lemma 4.1 depends on $N$ and does not correspond to a substitution in $t$.

## 4.3 The powers of $\pi$

**Lemma 4.5.** Let $k$ be the number of even entries of $q$, zero counted as even. Then the moment of Corollary 4.2 lies in $\pi^{-1}\mathbb{Q}$ if $k \in \{0,1\}$ and in $\mathbb{Q}$ if $k \in \{2,3\}$.

*Proof.* For $q_i$ even the argument $q_i + 1/2$ is a half-integer and $\Gamma$ contributes a rational multiple of $\sqrt{\pi}$; for $q_i$ odd it is a positive integer and $\Gamma$ contributes an integer. The numerator therefore carries $\pi^{k/2}$ times a rational. The number of odd entries is $3 - k$, so $Q \equiv 3 - k \pmod 2$; if $k$ is odd then $Q$ is even, $3 + Q/2$ is a half-integer and the denominator carries $\pi^{3/2}$, while if $k$ is even then $Q$ is odd and the denominator carries $\pi$. The quotient is $\pi^{(k-3)/2}$ for odd $k$ and $\pi^{(k-2)/2}$ for even $k$, which is $\pi^{-1}$ for $k \in \{0,1\}$ and $\pi^0$ for $k \in \{2,3\}$. ▫

**Corollary 4.6.** For every $N \geq 1$ the moment $m_N$ is a symmetric polynomial in $a_1, a_2, a_3$, homogeneous of degree $N$, whose coefficients are polynomials of degree at most one in $\pi^{-1}$ with rational coefficients. Moreover $m_1$ has degree zero in $\pi^{-1}$.

*Proof.* Expanding $w^N$ by the multinomial theorem and applying Corollary 4.2 term by term expresses $m_N$ as a rational combination of monomials $\prod_i a_i^{q_i}$ with $\sum_i q_i = N$, each weighted by a moment covered by Lemma 4.5. For $N = 1$ the only exponent triples are the permutations of $(1,0,0)$, for which $k = 2$, so no negative power of $\pi$ arises. ▫

**Remark 4.7.** The same conclusion follows from Proposition 4.3, by a different mechanism. Each factor $\mathrm{erf}(ta_i/\sqrt{2})$ equals $\pi^{-1/2}$ times an odd entire function of $t$ whose coefficients are free of $\pi$, while $\exp(1/2\, t^2 d^2)$ is even in $t$ and free of $\pi$. Expanding the product, a term involving exactly $k$ of the three error functions carries $\pi^{-k/2}$ and is a series in $t$ of parity $k$. Hence $\mathbb{E}[W^N]$ receives contributions only from $k \equiv N \pmod 2$, that is from $k \in \{0,2\}$ when $N$ is even and from $k \in \{1,3\}$ when $N$ is odd. Since $\mathbb{E}[|g|^N]$ is free of $\pi$ for even $N$ and is $\pi^{-1/2}$ times a $\pi$-free quantity for odd $N$, the quotient carries only $\pi^0$ and $\pi^{-1}$ in either case. For $N = 1$ the value $k = 3$ would require at least three powers of $t$, so only $k = 1$ occurs and the $\pi^{-1}$ term is absent.

We write $B_N$ for the coefficient of $\pi^{-1}$ in $m_N$, so that $B_1 = 0$.

## 4.4 The first five moments

Corollary 4.2 gives every moment as a finite sum. Expanding $w^N$ by the multinomial theorem and applying it term by term,

$$m_N = \frac{N!}{2\pi\,\Gamma(N+3/2)} \sum_{\substack{q_1+q_2+q_3=N \\ q_i \geq 0}} \prod_{i=1}^{3} \frac{a_i^{q_i}\,\Gamma(q_i+1/2)}{q_i!} \tag{6}$$

a closed expression from which each of the following may be obtained by hand. For $N = 1$ the three terms give $\Gamma(1)\Gamma(1/2)^2 a_i = \pi a_i$ and the prefactor is $1/(2\pi)$, returning $m_1 = L/2$.

**Proposition 4.8.** With $e_1, e_2, e_3$ the elementary symmetric functions of $a_1, a_2, a_3$,

$$m_1 = \frac{e_1}{2}, \qquad m_2 = \frac{e_1^2 - 2e_2}{3} + \frac{4e_2}{3\pi}, \qquad m_3 = \frac{2e_1^3 - 3e_1e_2 - 3e_3}{8} + \frac{3e_3}{2\pi}$$

$$m_4 = \frac{e_1^4 - 4e_1^2e_2 + 4e_2^2}{5} + \frac{8(2e_1^2e_2 - 4e_2^2 + e_1e_3)}{15\pi}$$

$$m_5 = \frac{8e_1^5 - 25e_1^3e_2 + 15e_1e_2^2 - 15e_1^2e_3 + 45e_2e_3}{48} + \frac{5e_3(e_1^2 - 2e_2)}{3\pi}$$

*Proof.* Formula (6) is a finite multinomial sum. For each $N \leq 5$, group its terms by the unordered exponent type $(q_1, q_2, q_3)$, evaluate the required gamma values, and reduce the resulting symmetric monomial sums to $e_1, e_2, e_3$ by Newton's identities. The complete hand calculation, including every exponent type and coefficient, is given in Appendix A. ▫

The first three of these formulas appear in [16] in projected-area coordinates; Proposition 4.3 gives every moment.

**Corollary 4.9.** Under the normalisation $e_1 = 1$, so that $e_2 = 1/2\,(1 - G)$ and $e_3 = J$,

$$B_2 = 2/3\,(1 - G), \qquad B_3 = 3/2\,J$$

In particular $B_2 > 0$ for every box, since $e_2 > 0$.

**Remark 4.10.** In the notation of Section 2.2 the first two moments read

$$\mathbb{E}[w] = \frac{L}{2}, \qquad \mathbb{E}[w^2] = \frac{L^2 - S}{3} + \frac{2S}{3\pi}$$

so that the first depends only on $V_1$ and the second only on $V_1$ and $V_2$. For the unit cube these give $3/2$ and $1 + 4/\pi$, the values obtained in [4].

**Remark 4.11.** Proposition 3.11 proves the total mass exactly through $\mathcal{P}(d) - \mathcal{P}(a_3) = \pi d/2$. For the unit cube, $d = \sqrt{3}$ and $d_i = \sqrt{2}$. As an independent check, numerical integration of Theorem 3.1 gives $\mathbb{E}[w] = 1.499999999999$ and $\mathbb{E}[w^2] = 2.273239544734$, against $3/2$ and $1 + 4/\pi = 2.273239544735$ in [4], an agreement to $10^{-12}$. The same density gives $\mathbb{E}[w^3] = 3.477464829$, $\mathbb{E}[w^4] = 5.365070725$ and $\mathbb{E}[w^5] = 8.341549431$, and the plateau of Proposition 3.3(1) is reproduced as $f \equiv 4/\sqrt{3} = 2.3094010768$ on $\left(\sqrt{2}, \sqrt{3}\right)$.

# 5. Structure of the moments

Section 4 gives the moments $m_N = \mathbb{E}[w^N]$ in closed form. Two of its results govern the present section. By Corollary 4.6, each $m_N$ is a symmetric polynomial in the edge lengths, homogeneous of degree $N$, whose coefficients are polynomials of degree at most one in $\pi^{-1}$ with rational coefficients, and $m_1 = 1/2 \sum_i a_i$ has degree zero; this rests on the parity count of Lemma 4.5.

The cumulants $\kappa_n$ of $w$ are defined by

$$\log \mathbb{E}[e^{tw}] = \sum_{n \geq 1} \kappa_n \frac{t^n}{n!}$$

and $g_n = \kappa_n / \kappa_1^n$ denotes the mean-normalised cumulants, which are invariant under $a \mapsto ca$. Throughout Sections 4 and 5 the edge lengths are treated as indeterminates and the moments and cumulants as elements of $\mathbb{Q}[a_1, a_2, a_3][\pi^{-1}]$; for such an element $F$ we write $[\pi^{-r}]F$ for the coefficient of $\pi^{-r}$ and $\deg_{\pi^{-1}} F$ for its degree as a polynomial in $\pi^{-1}$. With this convention $B_N = [\pi^{-1}] m_N$ is unambiguous, and $B_1 = 0$.

Throughout this section we normalise $\sum_i a_i = 1$, which is permitted by homogeneity; then $G = \sum_i a_i^2$, $J = a_1 a_2 a_3$ and $\kappa_1 = 1/2$, and by Corollary 4.9

$$B_2 = 2/3\,(1 - G), \qquad B_3 = 3/2\,J$$

## 5.1 The degree in $\pi^{-1}$

**Theorem 5.1.** For every $n \geq 2$,

$$\deg_{\pi^{-1}} \kappa_n = \lfloor \frac{n}{2} \rfloor \tag{7}$$

and the coefficient at top order is

$$\kappa_{2r}|_{\pi^{-r}} = (-1)^{r-1}(r-1)!\,\frac{(2r)!}{2^r r!}\,B_2^r$$

$$\kappa_{2r+1}|_{\pi^{-r}} = (-1)^{r-1}\binom{2r+1}{3}\frac{(2r-2)!}{2^{r-1}}\,B_2^{r-1}(B_3 - 3B_2 m_1)$$

*Proof.* In the moment-cumulant relation

$$\kappa_n = \sum_{\varpi} (-1)^{|\varpi|-1}\,(|\varpi|-1)! \prod_{B \in \varpi} m_{|B|}$$

summed over set partitions $\varpi$ of $\{1, \dots, n\}$, a factor $m_{|B|}$ can contribute at most $\pi^{-1}$ when $|B| \geq 2$, and contributes $\pi^0$ when $|B| = 1$, by Corollary 4.6. The degree of a term is therefore at most the number of blocks of size at least two, which for a partition of an $n$-element set is at most $\lfloor n/2 \rfloor$.

Suppose $n = 2r$. A partition with $r$ blocks of size at least two must have all blocks of size exactly two; there are $(2r)!/2^r r!$ such partitions, each with $|\varpi| = r$, and the first display follows.

Suppose $n = 2r + 1$. Blocks of size at least two, $r$ of them, occupy at least $2r$ elements, so exactly one of the following holds: there are $r - 1$ blocks of size two and one of size three, with no singleton; or there are $r$ blocks of size two and one singleton. The first type numbers $\binom{2r+1}{3}(2r-2)!/2^{r-1}(r-1)!$ and has $|\varpi| = r$; the second numbers $(2r+1)(2r)!/2^r r!$ and has $|\varpi| = r + 1$. Their contributions to $\kappa_n$ are therefore

$$c_a = (-1)^{r-1}\binom{2r+1}{3}\frac{(2r-2)!}{2^{r-1}}, \qquad c_b = (-1)^r (2r+1)\frac{(2r)!}{2^r}$$

and since $\binom{2r+1}{3} = (2r+1)(2r)(2r-1)/6$,

$$\frac{c_b}{c_a} = -\frac{(2r+1)(2r)(2r-1)}{2\binom{2r+1}{3}} = -3$$

independently of $r$, which gives the stated factorisation.

It remains to show that the top-order coefficient does not vanish. Since $1 - G = S/L^2 > 0$ we have $B_2 > 0$, so the even case is immediate. In the odd case

$$B_3 - 3B_2 m_1 = 3/2\,J - 3 \cdot 2/3\,(1-G) \cdot 1/2 = 3/2\,[J - 2/3\,(1-G)]$$

which vanishes only if $J = 2/3\,(1 - G)$, that is $V = 2/3\,SL$. Newton's inequality $S^2 \geq 12LV$ would then give $S \geq 8L^2$, contradicting $S \leq 2/3\,L^2$ from Lemma 2.2. Hence the degree is exactly $\lfloor n/2 \rfloor$ at every order. ▫

**Corollary 5.2.** No cumulant of order at most five contains $\pi^{-3}$, and the sixth does:

$$\kappa_6|_{\pi^{-3}} = 30\,B_2^3 = \frac{80}{9}(1 - G)^3, \qquad g_6|_{\pi^{-3}} = \frac{5120}{9}(1 - G)^3$$

**Corollary 5.3.** Taking $r = 2$,

$$\kappa_4|_{\pi^{-2}} = -3B_2^2, \qquad \kappa_5|_{\pi^{-2}} = -10\,B_2(B_3 - 3B_2 m_1)$$

whence, dividing by $\kappa_1^4$ and $\kappa_1^5$,

$$g_4|_{\pi^{-2}} = -\frac{64}{3}(1 - G)^2, \qquad g_5|_{\pi^{-2}} = \frac{320}{3}(1 - G)[2(1 - G) - 3J]$$

**Remark 5.4.** By Lemma 5.5 below, $g_4$ has four coefficients and $g_5$ has five. Corollary 5.3 determines the entire $\pi^{-2}$ component of each, and Corollary 5.2 gives the leading behaviour of $g_6$, in each case without computing the cumulant concerned. The argument uses only the parity count of Lemma 4.5, the combinatorics of set partitions, and Newton's inequalities.

## 5.2 Closure

**Lemma 5.5.** For every $n \geq 2$ the mean-normalised cumulant $g_n$ is a polynomial in $G$ and $J$ with coefficients in $\mathbb{Q}[\pi^{-1}]$, and the monomials that occur are among

$$\{\, G^\beta J^\gamma : \, 2\beta + 3\gamma \leq n \,\}$$

*Proof.* Each $m_N$ is symmetric in $a_1, a_2, a_3$, homogeneous of degree $N$, with coefficients in $\mathbb{Q}[\pi^{-1}]$; and $\kappa_n$ is a polynomial in $m_1, \dots, m_n$ with rational coefficients in which every term has total degree $n$. Hence $\kappa_n$ is symmetric and homogeneous of degree $n$ with coefficients in $\mathbb{Q}[\pi^{-1}]$.

The polynomials $p_1 = \sum_i a_i$, $p_2 = \sum_i a_i^2$ and $e_3 = a_1 a_2 a_3$ generate the symmetric polynomials in three variables over $\mathbb{Q}$, since $e_2 = (p_1^2 - p_2)/2$; their weights are 1, 2 and 3. Therefore

$$\kappa_n = \sum_{\alpha + 2\beta + 3\gamma = n} c_{\alpha\beta\gamma}\; p_1^\alpha p_2^\beta e_3^\gamma, \qquad c_{\alpha\beta\gamma} \in \mathbb{Q}[\pi^{-1}],\;\; \alpha, \beta, \gamma \geq 0$$

Dividing by $\kappa_1^n = (p_1/2)^n$ and substituting $\alpha = n - 2\beta - 3\gamma$,

$$\frac{p_1^\alpha p_2^\beta e_3^\gamma}{p_1^n} = \left(\frac{p_2}{p_1^2}\right)^\beta \left(\frac{e_3}{p_1^3}\right)^\gamma = G^\beta J^\gamma$$

and the requirement $\alpha \geq 0$ is exactly $2\beta + 3\gamma \leq n$. ▫

The argument is elementary and the result is stated as a lemma; its use is that the admissible monomials are known before any coefficient is computed. For $n = 2,3,4,5$ they are

$$\{1,G\},\qquad \{1,G,J\},\qquad \{1,G,G^2,J\},\qquad \{1,G,G^2,J,GJ\}$$

## 5.3 Closed forms

**Proposition 5.6.** Let $e_1 = L$, $e_2 = S/2$ and $e_3 = V$. Then

$$\kappa_2 = \frac{e_1^2}{12} - \frac{2e_2}{3} + \frac{4e_2}{3\pi}$$

$$\kappa_3 = \frac{5e_1e_2}{8} - \frac{2e_1e_2}{\pi} - \frac{3e_3}{8} + \frac{3e_3}{2\pi}$$

$$\kappa_4 = -\frac{e_1^4}{120} - \frac{43e_1^2e_2}{60} + \frac{12e_1^2e_2}{5\pi} + \frac{3e_1e_3}{4} - \frac{37e_1e_3}{15\pi} - \frac{8e_2^2}{15} + \frac{16e_2^2}{5\pi} - \frac{16e_2^2}{3\pi^2}$$

$$\kappa_5 = \frac{41e_1^3e_2}{48} - \frac{8e_1^3e_2}{3\pi} - \frac{15e_1^2e_3}{16} + \frac{17e_1^2e_3}{6\pi} + \frac{119e_1e_2^2}{48} - \frac{49e_1e_2^2}{3\pi} + \frac{80e_1e_2^2}{3\pi^2} - \frac{25e_2e_3}{16} + \frac{35e_2e_3}{3\pi} - \frac{20e_2e_3}{\pi^2}$$

and consequently

$$g_2 = \frac{4G-3}{3} + \frac{8(1-G)}{3\pi}$$

$$g_3 = -\frac{5G+6J-5}{2} + \frac{4(2G+3J-2)}{\pi}$$

$$g_4 = -\frac{2(16G^2-75G-90J+60)}{15} + \frac{16(12G^2-42G-37J+30)}{15\pi} - \frac{64(1-G)^2}{3\pi^2}$$

$$g_5 = \frac{119G^2+150GJ-320G-330J+201}{6} - \frac{8(49G^2+70GJ-114G-104J+65)}{3\pi} + \frac{320(G-1)(2G+3J-2)}{3\pi^2}$$

*Proof.* Insert the raw moments of Proposition 4.8 into the standard moment–cumulant identities through order five and collect separately the coefficients of each power of $\pi^{-1}$. Every resulting expression is symmetric and homogeneous in $a_1, a_2, a_3$, and hence reduces uniquely to a polynomial in $e_1, e_2, e_3$. Appendix A carries out this hand reduction coefficient by coefficient. Finally substitute

$$e_2 = 1/2\, e_1^2(1-G),\qquad e_3 = Je_1^3,\qquad \kappa_1 = 1/2\, e_1,$$

and divide by $\kappa_1^n$. The powers of $e_1$ cancel, giving the displayed formulas for $g_2, g_3, g_4, g_5$. ▫

**Remark 5.7.** The $\pi^{-2}$ components of $g_4$ and $g_5$ agree with Corollary 5.3, and no coefficient of $g_5$ carries $\pi^{-3}$, as Corollary 5.2 requires. Since Proposition 5.6 is obtained independently of Theorem 5.1, this is a check on both.

### 5.4 Consequences for the variance

**Corollary 5.8.** In terms of the intrinsic volumes,

$$\mathbb{E}[w] = \frac{V_1}{2}, \qquad \mathrm{Var}(w) = \frac{V_1^2}{12} - \left(\frac{1}{3} - \frac{2}{3\pi}\right) 2V_2$$

so the mean depends only on $V_1$ and the variance only on $V_1$ and $V_2$, whereas $\kappa_3$ involves $V_3$ with coefficient $3/2\,\pi^{-1} - 3/8 \neq 0$.

The variance formula appears in [16] in projected-area coordinates.

**Corollary 5.9.** The squared coefficient of variation of $w$ is

$$\mathrm{CV}^2 = g_2 = \frac{4(\pi - 2)G + (8 - 3\pi)}{3\pi}$$

which is strictly increasing in $G$; hence, by Section 2.3,

$$\frac{16 - 5\pi}{9\pi} \leq \mathrm{CV}^2 < \frac{1}{3}$$

the lower bound attained precisely for the cube. In particular no box has width of constant length, and among boxes the cube has the least relative dispersion of width.

*Proof.* The coefficient of $G$ is $4(\pi - 2)/3\pi > 0$. Substituting $G = 1/3$ gives $(16 - 5\pi)/9\pi = 0.0103287$, attained only at the cube by Section 2.3, and $G \to 1$ gives $1/3$, which is not attained. ▫

**Remark 5.10.** Since $G = (d/L)^2$ by Section 2.3, the squared coefficient of variation is a strictly increasing function of the ratio of the space diagonal to the sum of the edge lengths.

## 6. The graded dependence on the intrinsic volumes

Corollary 5.8 exhibits a graded dependence: the mean of $w$ is a function of $V_1$ alone, the variance of $V_1$ and $V_2$, and the third cumulant involves $V_3$. This section shows that the grading is sharp, so that no step may be omitted, and that three invariants suffice.

### 6.1 The admissible invariants

**Lemma 6.1.** Fix $L > 0$ and $S > 0$, and put

$$g(z) = z^3 - Lz^2 + S/2\,z, \qquad z_\pm = \frac{L \pm \sqrt{L^2 - 3/2\,S}}{3}$$

A box with these values of $V_1$ and $V_2$ exists if and only if $S \leq 2/3\,L^2$, and the admissible volumes are then

$$0 < V \leq g(z_-) \quad \text{if } 0 < S \leq 1/2\,L^2, \qquad\qquad g(z_+) \leq V \leq g(z_-) \quad \text{if } 1/2\,L^2 < S \leq 2/3\,L^2$$

The value $V = 0$ is never attained. The upper endpoint corresponds to a box with a repeated edge, as does the lower endpoint when it is positive, and at $S = 2/3\,L^2$ the two coincide and give the cube.

*Proof.* By Lemma 2.1 the edges are the roots of $p(z) = g(z) - V$, so a box with the given invariants exists precisely when $p$ has three positive roots. Since $p'(z) = 3z^2 - 2Lz + S/2$, real critical points require $L^2 \geq 3/2\,S$, which is Newton's first inequality of Lemma 2.2. Given that, $p$ has three real roots if and only if $g(z_+) \leq V \leq g(z_-)$, with equality at either end giving a double root. Moreover, if $V > 0$ then $p(z) < 0$ for every $z \leq 0$, since $z^3$, $-Lz^2$ and $S/2\,z$ are all non-positive there while $-V < 0$; hence $p$ has no non-positive root, and any three real roots are all positive. Since a box requires $V > 0$, this is the only further condition.

It remains to locate $g(z_+)$. At a critical point $S/2 = 2Lz - 3z^2$, whence

$$g(z) = z^3 - Lz^2 + (2Lz - 3z^2)z = z^2(L - 2z)$$

and $L - 2z_+ = 1/3\left(L - 2\sqrt{L^2 - 3/2\,S}\right)$, which is non-positive exactly when $L^2 \leq 4L^2 - 6S$, that is when $S \leq 1/2\,L^2$. In that range the constraint $V \geq g(z_+)$ is vacuous and $V > 0$ binds instead. ▫

**Remark 6.2.** In the scale-invariant variables of Section 2.3 the threshold $S = 1/2\,L^2$ is $G = 1/2$. The closure of the admissible region in the $(G, J)$ plane therefore contains the segment $J = 0$, $1/2 \leq G \leq 1$, along which the volume degenerates: as a single edge tends to zero the limit point traverses $1/2 \leq G < 1$, and $G = 1$ is reached when two edges do. The curved part of the boundary is the two-equal-edges locus.

## 6.2 The second step is sharp

**Proposition 6.3.** Fix $L$. Then $\mathrm{Var}(w)$ is a strictly decreasing function of $S$, and as $S$ ranges over $(0, 2/3\,L^2]$ it takes every value in

$$[\, L^2(1/12 - 2/9 + 4/9\pi),\ \ L^2/12\,) = [\,0.0025821716\,L^2,\ 0.0833333333\,L^2\,)$$

the lower endpoint attained only by the cube. In particular $V_1$ does not determine $\mathrm{Var}(w)$.

*Proof.* By Corollary 5.8, $\mathrm{Var}(w) = L^2/12 - cS$ with $c = 1/3 - 2/3\pi = 0.1211267425 > 0$, so the dependence on $S$ is strictly decreasing. Newton's first inequality gives $S \leq 2/3\,L^2$ with equality only when the edges coincide, and substituting the endpoints gives the stated interval. The unattained upper endpoint is approached, for example, by the boxes with edges $(L - 2\varepsilon, \varepsilon, \varepsilon)$ as $\varepsilon \downarrow 0$. ▫

**Example 6.4.** The boxes $(1,1,1)$ and $(1/2, 1, 3/2)$ both have $L = 3$, while $S = 6$ and $S = 11/2$ respectively, so that

$$\mathrm{Var}(w) = 0.0232395447 \quad \text{and} \quad 0.0838029160$$

## 6.3 The third step is sharp

**Proposition 6.5.** Fix $L$ and $S$ with $S < 2/3\,L^2$. Then $g_3$ is a strictly increasing function of $V$ on the interval of Lemma 6.1, and in particular $(V_1, V_2)$ does not determine the third cumulant.

*Proof.* By Proposition 5.6, $g_3$ depends on $V$ only through $J = V/L^3$, with coefficient $3(4/\pi - 1) = 0.8197186$, which is positive. The interval of Lemma 6.1 is non-degenerate because $S < 2/3\,L^2$ makes $z_+ \neq z_-$. ▫

**Example 6.6.** Take $L = 6$ and $S = 22$, so that $G = 1 - S/L^2 = 7/18$. By Lemma 6.1 the admissible volumes are

$$V \in \left[\, 6 - 2\sqrt{3}/9\,,\ 6 + 2\sqrt{3}/9 \,\right] = [5.6150998205,\ 6.3849001795]$$

and the two endpoints are attained by the boxes with a repeated edge, labelled by the sign in $V$, so that $B_+$ carries the larger volume and hence the double root $z_-$,

$$B_+ = \left(2 - \sqrt{3}/3\,,\ 2 - \sqrt{3}/3\,,\ 2 + 2\sqrt{3}/3\right), \qquad B_- = \left(2 + \sqrt{3}/3\,,\ 2 + \sqrt{3}/3\,,\ 2 - 2\sqrt{3}/3\right)$$

whose edges are, to ten decimal places,

$$(1.4226497308,\ 1.4226497308,\ 3.1547005384)$$

$$(2.5773502692,\ 2.5773502692,\ 0.8452994616)$$

Both have $L = 6$ and $S = 22$, hence the same mean and the same variance of width, while

$$g_3(B_+) = -0.0041732322, \qquad g_3(B_-) = -0.0070946197$$

the difference being exactly $\sqrt{3}\,(4 - \pi)/162\pi = 0.0029213875$.

Figure 2 shows the admissible region and both witness pairs.

**Figure 2.** The admissible region in the $(G, J)$ plane, with $G = 1 - S/L^2$ and $J = V/L^3$. Its curved boundary is the image of the two-equal-edges locus $a = (1,1,t)$. The two boundary arcs meet in a semicubical cusp at the cube $(1/3, 1/27)$, indicated by the star, and their common tangent there has slope $-1/6$. By Remark 6.2 the closure also contains the segment $J = 0$, $1/2 \leq G \leq 1$. Newton's inequality $J \leq (1 - G)^2/12$ of Lemma 2.2 is dashed. Among nondegenerate boxes equality occurs only at the cube; the Newton curve also meets the degenerate closure at $(G, J) = (1,0)$ and is otherwise strictly weaker. The star and the open circle also represent the boxes of Example 6.4, which share $L$ and differ in $S$: their horizontal separation is the variance step of Proposition 6.3. The squares are $B_\pm$ of Example 6.6, which share $L$ and $S$ and differ only in $V$; the segment joining them is the entire admissible range of $V$ at $G = 7/18$, and is the third-cumulant step of Proposition 6.5. The open circle lies on that segment.

## 6.4 Three invariants suffice

**Proposition 6.7.** The triple $(V_1, V_2, V_3)$ determines the law of $w$ completely.

*Proof.* By Lemma 2.1 the triple determines $a_1, a_2, a_3$ up to permutation, and $w$ is symmetric in the edges. ▫

## 6.5 The first three cumulants determine the law

**Theorem 6.8.** The first three cumulants of $w$ determine the edge lengths up to permutation, and hence the whole distribution. Explicitly,

$$L = 2\kappa_1, \qquad S = \frac{3\pi}{\pi - 2}\left(\frac{L^2}{12} - \kappa_2\right), \qquad V = \frac{8\pi}{3(4-\pi)}\left(\kappa_3 - \frac{5\pi - 16}{16\pi} LS\right)$$

after which the edges are the positive roots of $z^3 - Lz^2 + S/2\, z - V$.

*Proof.* By Corollary 5.8, $\kappa_1 = L/2$. Proposition 5.6 in the intrinsic volumes reads

$$\kappa_2 = \frac{L^2}{12} - \frac{\pi - 2}{3\pi} S, \qquad \kappa_3 = \frac{5\pi - 16}{16\pi} LS + \frac{3(4 - \pi)}{8\pi} V$$

Since $\pi - 2 > 0$ the second determines $S$ from $L$ and $\kappa_2$, and since $4 - \pi > 0$ the third determines $V$ from $L$, $S$ and $\kappa_3$. Lemma 2.1 then recovers the edges. ▫

**Corollary 6.9.** A triple $(\kappa_1, \kappa_2, \kappa_3)$ of real numbers is the triple of first three cumulants of the width of a rectangular box if and only if the quantities

$$L = 2\kappa_1, \qquad S = \frac{3\pi}{\pi - 2}\left(\frac{L^2}{12} - \kappa_2\right), \qquad V = \frac{8\pi}{3(4 - \pi)}\left(\kappa_3 - \frac{5\pi - 16}{16\pi} LS\right)$$

satisfy $L > 0$, $0 < S \leq 2/3\, L^2$, and

$$0 < V \leq g(z_-) \quad \text{if } 0 < S \leq 1/2\, L^2, \qquad\qquad g(z_+) \leq V \leq g(z_-) \quad \text{if } 1/2\, L^2 < S \leq 2/3\, L^2$$

in the notation of Lemma 6.1. The box is then determined up to congruence.

*Proof.* The three displayed expressions invert the relations of Theorem 6.8, so a box with the given cumulants exists precisely when the recovered $L$, $S$ and $V$ are the intrinsic volumes of a box. By Lemma 2.1 that holds precisely when $z^3 - Lz^2 + S/2\, z - V$ has three positive roots, and Lemma 6.1 states exactly when it does. Uniqueness is Theorem 6.8. ▫

**Corollary 6.10.** For every $n \geq 4$ there is a polynomial $P_n$ over $\mathbb{Q}(\pi)$ with $\kappa_n = P_n(\kappa_1, \kappa_2, \kappa_3)$, and similarly for the raw moments. Propositions 6.3 and 6.5 show that $\kappa_1$ alone does not determine $\kappa_2$ and that $(\kappa_1, \kappa_2)$ does not determine $\kappa_3$, so the initial cumulant sequence $(\kappa_1, \kappa_2, \kappa_3)$ is minimal in length. This is a statement about initial segments, and does not assert that no other pair of statistics determines the box.

*Proof.* By the proof of Lemma 5.5 every $\kappa_n$ is a polynomial in $L$, $S$ and $V$ over $\mathbb{Q}(\pi)$, and by Theorem 6.8 these are polynomials in $\kappa_1, \kappa_2, \kappa_3$ over the same field. ▫

**Remark 6.11.** Propositions 6.3, 6.5 and 6.7 with Theorem 6.8 locate the width distribution exactly within the intrinsic volumes: the nested hierarchy $V_1$, $(V_1, V_2)$, $(V_1, V_2, V_3)$ is strict and closes at its third step, and no additional invariant is needed to determine the width law. Lemma 5.5 records the corresponding algebraic statement, that every mean-normalised cumulant is a polynomial in the two scale-invariant combinations $G$ and $J$. Section 7 shows that the corresponding closure fails for general convex bodies, for the width and for the brightness alike.

# 7. Beyond the box

## 7.1 Moments of the brightness for a general convex body

This section carries out two separate extensions. Sections 7.1 and 7.2 compute the low moments of the brightness of an arbitrary convex body and show that the closure of Section 6 in the intrinsic volumes does not survive. Section 7.3 characterises exactly when the sign-pattern representation underlying the box calculation survives. The cellwise co-area method itself applies more broadly to every full-dimensional polytope, as explained in Section 8.

By Section 1.3 the brightness of a convex body is the cosine transform of its surface area measure, $A_K = 1/2\ CS_K$ where $(C\mu)(u) = \int |u \cdot n|\, d\mu(n)$. The first three moments are therefore available in general.

The bivariate Gaussian absolute-moment identity underlying the following calculation is classical; see Kamat [8] and Nabeya [10]. We include a short first-principles derivation because its angular normalisation is used repeatedly below.

**Theorem 7.1.** For $u$ uniform on $S^2$ and unit vectors $n, n'$ at angle $\theta$,

$$\langle |u \cdot n|\, |u \cdot n'| \rangle = \Omega(\theta) := \frac{2}{3\pi}\left[\sin\theta + \cos\theta \left(\frac{\pi}{2} - \theta\right)\right]$$

*Proof.* With $g$ standard normal in $\mathbb{R}^3$ and $u = g/|g|$ as in Lemma 4.1, the variables $X = g \cdot n$ and $Y = g \cdot n'$ are jointly standard normal with correlation $\cos\theta$, and $|g|$ is independent of $u$, so $\mathbb{E}|XY| = \mathbb{E}|g|^2 \cdot \Omega = 3\Omega$.

To evaluate $\mathbb{E}|XY|$, realise the pair as $X = Z_1$ and $Y = \cos\theta\, Z_1 + \sin\theta\, Z_2$ with $(Z_1, Z_2)$ standard bivariate normal. In polar form $Z_1 = R\cos\varphi$ and $Z_2 = R\sin\varphi$ with $\varphi$ uniform on $[0,2\pi)$ independent of $R$, and $\mathbb{E}[R^2] = 2$; then $Y = R\cos(\varphi - \theta)$ and

$$\mathbb{E}|XY| = \mathbb{E}[R^2] \cdot \frac{1}{2\pi}\int_0^{2\pi} |\cos\varphi \cos(\varphi - \theta)|\, d\varphi = \frac{1}{\pi}\int_0^{2\pi} |\cos\varphi \cos(\varphi - \theta)|\, d\varphi$$

The integrand has period $\pi$, so the integral over $[0,2\pi)$ is twice that over $[0,\pi]$. On $[0,\pi]$ the factor $\cos\varphi$ changes sign only at $\pi/2$ and $\cos(\varphi - \theta)$ only at $\theta + \pi/2$; taking $\theta \in [0, \pi/2]$, so that both lie in the range, the product is non-negative on $[0, \pi/2]$ and on $[\theta + \pi/2\,, \pi]$ and non-positive between them. From $2\cos\varphi\cos(\varphi - \theta) = \cos\theta + \cos(2\varphi - \theta)$ the antiderivative is

$$P(\varphi) = \frac{\varphi\cos\theta}{2} + \frac{\sin(2\varphi - \theta)}{4}$$

so the three pieces contribute $P(\pi/2) - P(0)$, then $-\big(P(\theta + \pi/2) - P(\pi/2)\big)$, then $P(\pi) - P(\theta + \pi/2)$. The four values are

$$P(0) = -\sin\theta/4\,, \quad P(\pi/2) = \pi\cos\theta/4 + \sin\theta/4\,, \quad P(\theta + \pi/2) = (2\theta + \pi)\cos\theta/4 - \sin\theta/4\,, \quad P(\pi) = \pi\cos\theta/2 - \sin\theta/4$$

and their alternating sum, doubled, is $(\pi - 2\theta)\cos\theta + 2\sin\theta$. Hence

$$\mathbb{E}|XY| = \frac{1}{\pi}[(\pi - 2\theta)\cos\theta + 2\sin\theta] = \frac{2}{\pi}\left[\sin\theta + \cos\theta\left(\frac{\pi}{2} - \theta\right)\right]$$

If $\theta \in [\pi/2\,, \pi]$, replace $n'$ by $-n'$. The absolute product $|u \cdot n|\, |u \cdot n'|$ is unchanged, while the angle between $n$ and $-n'$ is $\pi - \theta \in [0, \pi/2]$. Applying the formula already proved and then substituting $\pi - \theta$ gives

$$\sin(\pi - \theta) + \cos(\pi - \theta)\left(\frac{\pi}{2} - (\pi - \theta)\right) = \sin\theta + \cos\theta\left(\frac{\pi}{2} - \theta\right),$$

so the same expression holds on the whole interval $[0,\pi]$. Dividing by 3 gives $\Omega$. ▫

Two values are immediate: $\Omega(0) = \Omega(\pi) = 1/3$ and $\Omega(\pi/2) = 2/3\pi$, recovering $\langle u_i^2 \rangle$ and $\langle |u_i u_j| \rangle$. Moreover $\Omega$ is symmetric about $\pi/2$, since

$$\Omega(\pi - \theta) = 2/3\pi\, [\sin\theta + \cos\theta(\pi/2 - \theta)] = \Omega(\theta)$$

so a pair of normals meeting at $\pi - \alpha$ contributes as if meeting at $\alpha$.

The trivariate analogue is also available, through the classical Gaussian absolute product moment in three variables.

**Proposition 7.2.** Let $n_1, n_2, n_3$ be unit vectors and put $\rho_{ij} = n_i \cdot n_j$.

1. *Suppose they are pairwise non-parallel, and define $A_k \in [0,\pi]$ by*

$$\cos A_k = \frac{\cos\theta_{ij} - \cos\theta_{ik}\cos\theta_{jk}}{\sin\theta_{ik}\sin\theta_{jk}}$$

*where $\theta_{ij}$ is the angle between $n_i$ and $n_j$. When the $n_i$ are linearly independent these are the interior angles of the spherical triangle with vertices $n_1, n_2, n_3$, and when they are coplanar each $A_k$ is 0 or $\pi$. Then for $u$ uniform on $S^2$,*

$$\Omega_3(n_1,n_2,n_3) := \mathbb{E}\left[\prod_{i=1}^{3} |u \cdot n_i|\right] = \frac{1}{4\pi}\left[|\det(n_1,n_2,n_3)| + \sum_k \left(\rho_{ij} + \rho_{ik}\rho_{jk}\right)\left(\frac{\pi}{2} - A_k\right)\right]$$

*with $\{i,j,k\} = \{1,2,3\}$.*

2. *If two normals are parallel, say $n_1 = \pm n_2$, and $\theta$ is the angle between $n_1$ and $n_3$, then $\Omega_3 = 1/8\,(1 + \cos^2\theta)$.*

*Proof.* With $g$ standard normal in $\mathbb{R}^3$ and $u = g/|g|$ as in Lemma 4.1, the variables $X_i = g \cdot n_i$ are standard normal with correlations $\rho_{ij}$, and $|g|$ is independent of $u$, so $\mathbb{E}|X_1X_2X_3| = \mathbb{E}|g|^3\, \Omega_3 = 8\sqrt{2/\pi}\, \Omega_3$. Nabeya [10, p. 19] evaluates the left side as

$$(\frac{2}{\pi})^{3/2}\left[\sqrt{\det R} + \sum_k \left(\rho_{ij} + \rho_{ik}\rho_{jk}\right)\arcsin\rho_{ij\cdot k}\right]$$

with $R = \left(\rho_{ij}\right)$ and $\rho_{ij\cdot k}$ the partial correlation; dividing by $8\sqrt{2/\pi}$ produces the factor $1/(4\pi)$. Since $R$ is the Gram matrix of the $n_i$ we have $\det R = \det(n_1,n_2,n_3)^2$, and the partial correlation $\rho_{ij\cdot k}$ is by definition the right side of the displayed relation, so $\rho_{ij\cdot k} = \cos A_k$ and $\arcsin\rho_{ij\cdot k} = \pi/2 - A_k$. This proves the formula whenever $R$ is non-singular, that is whenever $n_1, n_2, n_3$ are linearly independent.

Suppose now that they are coplanar but still pairwise non-parallel. Choose linearly independent unit vectors $n_i^{(\varepsilon)} \to n_i$, and realise all the Gaussian vectors on one probability space by setting $X_i^{(\varepsilon)} = g \cdot n_i^{(\varepsilon)}$ with a single $g$. Then $X_i^{(\varepsilon)} \to X_i$ pointwise, and $\prod_i \left|X_i^{(\varepsilon)}\right| \le |g|^3$ for every $\varepsilon$, which is integrable; by dominated convergence the left sides converge. The right side is a continuous function of the Gram entries $\rho_{ij}$ on the region where no two are $\pm 1$, since each $A_k$ is continuous there. Letting $\varepsilon \to 0$ gives the formula in the coplanar case, where each $A_k$ is 0 or $\pi$. If $n_1 = \pm n_2$ the product is $(u \cdot n_1)^2 |u \cdot n_3|$; writing $X = g \cdot n_1$ and $Y = g \cdot n_3$ with correlation $\rho$ and $X = \rho Y + \sqrt{1-\rho^2}\, Z$ with $Z$ standard normal independent of $Y$,

$$\mathbb{E}[X^2 |Y|] = \rho^2\, \mathbb{E}|Y|^3 + (1-\rho^2)\mathbb{E}|Y| = \sqrt{2/\pi}\,(1+\rho^2)$$

and division by $8\sqrt{2/\pi}$ gives the stated value; this covers three parallel normals as the case $\cos^2\theta = 1$. ▫

**Corollary 7.3.** For every convex body $K$,

$$\mathbb{E}[A_K] = \frac{S(K)}{4}, \qquad \mathbb{E}[A_K^2] = \frac{1}{4}\iint_{S^2\times S^2} \Omega\,(\theta_{nn'})\, dS_K(n)\, dS_K(n')$$

$$\mathbb{E}[A_K^3] = \frac{1}{8}\iiint_{(S^2)^3} \Omega_3\,(n,n',n'')\, dS_K(n)\, dS_K(n')\, dS_K(n'')$$

For a ball of radius $R$ the three right sides are $\pi R^2$, $\pi^2 R^4$ and $\pi^3 R^6$, consistent with $A_K \equiv \pi R^2$. For a polytope the measure is atomic and the degenerate clause of Proposition 7.2 supplies the terms in which two or three of the normals coincide.

**Remark 7.4.** For a box, $A_K$ is by (1) the width of the box with edges $a_j a_k$, so Corollary 7.3 must reproduce $m_1, m_2$ and $m_3$ of Proposition 4.8. It does, the third identity testing the degenerate terms of Proposition 7.2 as well as the generic ones.

## 7.2 Closure fails outside the box

The closure of Section 6 in the intrinsic volumes is special to boxes. Already at second order, neither the width variance nor the brightness variance of a general convex body is determined by $V_1, V_2, V_3$.

**Theorem 7.5.** *Among convex bodies $K \subset \mathbb{R}^3$, the intrinsic volumes $V_1(K), V_2(K), V_3(K)$ determine neither* $\mathrm{Var}(w_K)$ *nor* $\mathrm{Var}(A_K)$.

The required pair of bodies is constructed explicitly next.

**Example 7.6.** A right prism of height $h$ over a triangle of area $A$ and perimeter $P$ has

$$V = Ah, \qquad S = 2A + Ph, \qquad V_1 = h + 1/2\, P$$

The last identity follows from the edge-angle formula for $V_1$. If the interior angles of the triangle are $\alpha_1, \alpha_2, \alpha_3$, the three vertical edges contribute

$$h\sum_{i=1}^{3}(\pi - \alpha_i) = 2\pi h$$

to $2\pi V_1$, while the six horizontal edges contribute $\pi P$. Hence

$$2\pi V_1 = 2\pi h + \pi P$$

and $V_1 = h + P/2$.

Take

$$h = 1, \qquad A = 4, \qquad P = 10$$

Then every such prism has

$$V = 4, \qquad S = 18, \qquad V_1 = 6$$

These are also the intrinsic volumes of the box with edges $(1,1,4)$.

We first compare the brightness variances. For the box,

$$A_K(u) = 4|u_1| + 4|u_2| + |u_3|$$

so Proposition 4.8, applied with coefficients $(4,4,1)$, gives

$$\mathbb{E}[A_K] = \frac{9}{2}, \qquad \mathrm{Var}(A_K) = \frac{32}{\pi} - \frac{37}{4}$$

For a right prism of height $h$ over a triangle with side lengths $L_m$ and interior angles $\alpha_k$, the two triangular faces have normals $\pm q$ and area $A$, while the three rectangular faces have in-plane normals and areas $hL_m$. The in-plane normals meet $\pm q$ at angle $\pi/2$, and the normals corresponding to sides $m, n$ meet at $\pi - \alpha_k$, where $\{m, n, k\} = \{1,2,3\}$. Since $\Omega(\pi - \alpha) = \Omega(\alpha)$, Corollary 7.3 gives

$$\mathrm{Var}(A_K) = \frac{1}{4}\left[\frac{4}{3}A^2 + \frac{8AhP}{3\pi} + h^2\left(\frac{1}{3}\sum_m L_m^2 + 2\sum_{m<n} L_m L_n \Omega(\alpha_k)\right)\right] - \left(\frac{S}{4}\right)^2$$

Choose the right triangle whose legs and hypotenuse are

$$a = \frac{29 + \sqrt{41}}{10}, \qquad b = \frac{29 - \sqrt{41}}{10}, \qquad c = \frac{21}{5}$$

Indeed,

$$a + b = \frac{29}{5}, \qquad ab = 8$$

so its perimeter is 10 and its area is $ab/2 = 4$; moreover

$$a^2 + b^2 = (a + b)^2 - 2ab = \frac{841}{25} - 16 = \frac{441}{25} = c^2$$

so it is right-angled. Let $\alpha$ be the angle opposite $a$ and put $\beta = \pi/2 - \alpha$. First,

$$\sum_m L_m^2 = a^2 + b^2 + c^2 = 2c^2 = \frac{882}{25}.$$

Moreover, $\sin\alpha = a/c$, $\cos\alpha = b/c$, $\sin\beta = b/c$ and $\cos\beta = a/c$, so the formula for $\Omega$ gives

$$\begin{aligned}\sum_{m<n} L_m L_n \Omega(\alpha_k) &= bc\,\Omega(\alpha) + ca\,\Omega(\beta) + ab\,\Omega(\pi/2) \\ &= \frac{2}{3\pi}(3ab + a^2\alpha + b^2\beta) \\ &= \frac{2}{3\pi}(24 + a^2\alpha + b^2\beta).\end{aligned}$$

Direct substitution into the brightness formula now gives

$$\begin{aligned}\mathrm{Var}(A_K) &= \frac{1}{4}\left[\frac{64}{3} + \frac{320}{3\pi} + \frac{294}{25} + \frac{4}{3\pi}(24 + a^2\alpha + b^2\beta)\right] - \frac{81}{4} \\ &= -\frac{3593}{300} + \frac{104 + a^2\alpha + b^2\beta}{3\pi}.\end{aligned}$$

Since $a > b$, one has $\alpha > \pi/4$. Also $a^2 < 3b^2$, so

$$\tan\alpha = \frac{a}{b} < \sqrt{3}$$

and therefore $\alpha < \pi/3$. Now

$$a^2\alpha + b^2\beta = \frac{\pi}{2}b^2 + (a^2 - b^2)\alpha$$

is increasing in $\alpha$, and

$$a^2 - b^2 = \frac{29\sqrt{41}}{25}$$

Using $\alpha < \pi/3$ yields

$$\mathrm{Var}(A_K) < \frac{104}{3\pi} - \frac{788}{75} + \frac{29\sqrt{41}}{900}$$

Subtracting this upper bound from the box value gives

$$\mathrm{Var}(A_{\text{box}}) - \mathrm{Var}\big(A_{\text{prism}}\big) > \frac{377}{300} - \frac{8}{3\pi} - \frac{29\sqrt{41}}{900}$$

The elementary bounds $\pi > 3$ and $\sqrt{41} < 13/2$ give

$$\frac{8}{3\pi} + \frac{29\sqrt{41}}{900} < \frac{8}{9} + \frac{377}{1800} = \frac{1977}{1800} < \frac{2262}{1800} = \frac{377}{300}$$

and hence

$$\mathrm{Var}(A_{\text{box}}) - \mathrm{Var}\big(A_{\text{prism}}\big) > \frac{19}{120} > 0$$

This proves that the brightness variances are different without numerical evaluation. For orientation, their numerical values are respectively 0.9359163579 and 0.6993376584.

We next compare the width variances for the same two bodies. For the box,

$$w_K(u) = |u_1| + |u_2| + 4|u_3|$$

and Proposition 4.8 gives

$$\mathbb{E}[w_K] = 3, \qquad \operatorname{Var}(w_K) = \frac{12}{\pi} - 3$$

Let $T$ be the triangular base of the prism, with cyclic edge vectors $E_1, E_2, E_3$. For $x_i = v_i \cdot u$, the elementary identity

$$\max_i x_i - \min_i x_i = \frac{1}{2} \sum_{\text{cyc}} |x_i - x_j|$$

gives

$$w_T(u) = \frac{1}{2} \sum_{m=1}^{3} |u \cdot E_m|$$

Since the prism is the Minkowski sum of $T$ and a segment of length $h$ parallel to $q$,

$$w_K(u) = h|u \cdot q| + \frac{1}{2} \sum_{m=1}^{3} |u \cdot E_m|$$

The four coefficient vectors have lengths $h, L_1/2, L_2/2, L_3/2$. The vector $hq$ is orthogonal to the other three, while the cyclic edge vectors $E_m, E_n$ meet at $\pi - \alpha_k$. Applying Theorem 7.1 term by term therefore gives

$$\mathbb{E}[w_K^2] = \frac{h^2}{3} + \frac{\sum_m L_m^2}{12} + \frac{2hP}{3\pi} + \frac{1}{2} \sum_{m<n} L_m L_n \Omega(\alpha_k)$$

For $h = 1$ and the same right triangle, $\mathbb{E}[w_K] = 1/2\,(h + P/2) = 3$. Using the two sums evaluated above gives

$$\begin{aligned} \mathbb{E}[w_K^2] \quad &= \frac{1}{3} + \frac{1}{12}\frac{882}{25} + \frac{20}{3\pi} + \frac{24 + a^2\alpha + b^2\beta}{3\pi} \\ &= \frac{491}{150} + \frac{44 + a^2\alpha + b^2\beta}{3\pi}, \end{aligned}$$

and therefore

$$\operatorname{Var}(w_K) = -\frac{859}{150} + \frac{44 + a^2\alpha + b^2\beta}{3\pi}.$$

Using the same bounds $\pi/4 < \alpha < \pi/3, \pi > 3$ and $\sqrt{41} < 13/2$ gives

$$\operatorname{Var}(w_{\text{box}}) - \operatorname{Var}(w_{\text{prism}}) > \frac{19}{120} > 0$$

Thus the width variances also differ. In fact the two variance gaps are equal, both being

$$\frac{409}{150} - \frac{8 + a^2\alpha + b^2\beta}{3\pi}$$

At height one the brightness and width of the prism have the same in-plane coefficients $L_m/2$ and differ only in the coefficient of $|u \cdot q|$, so their variances differ by a quantity determined by the common intrinsic volumes. This equality of the two gaps is particular to $h = 1$ and is not a general identity.

*Proof of Theorem 7.5.* The box and the prism constructed in Example 7.6 have the same three intrinsic volumes. The two exact inequalities established there show that their brightness variances differ and that their width variances differ. Consequently neither variance is determined by $V_1, V_2, V_3$. ▫

The mechanism is visible in the formulas. By Corollary 7.3 the brightness variance depends on the pairwise angles between the normals carrying $S_K$, while the width variance of a finite sum of absolute linear forms depends, through Theorem 7.1, on the pairwise angles between its coefficient vectors. Neither family of angular data is determined by the intrinsic volumes.

## 7.3 How far the density argument reaches

We now return from moment closure to the density argument, and identify precisely when its sign-pattern representation survives.

**Proposition 7.7.** Let $DK = 1/2\,(K - K)$ denote the central symmetral of $K$. For vectors $g_1, \dots, g_N$ the conditions

$$w_K(u) = \sum_{k=1}^{N} |\, u \cdot g_k| \qquad \text{and} \qquad DK = \sum_{k=1}^{N} [-1/2\, g_k,\ 1/2\, g_k]$$

are equivalent. Consequently $w_K$ is a finite sum of absolute values of linear forms if and only if $K - K$, equivalently $DK$, is a zonotope.

*Proof.* $h_{DK}(u) = 1/2\, h_{K-K}(u) = 1/2\,\big(h_K(u) + h_K(-u)\big) = 1/2\, w_K(u)$, while the support function of $\sum_k [-1/2\, g_k, 1/2\, g_k]$ is $1/2 \sum_k |u \cdot g_k|$. The equivalence follows from the uniqueness of support functions. ▫

A box is a zonotope with three mutually orthogonal generators, and $DK$ is then the box itself.

**Remark 7.8.** The class of bodies satisfying Proposition 7.7 is strictly larger than the zonotopes, and contains bodies that are not centrally symmetric. Let $K$ be the right prism of height $h$ over the triangle with vertices $v_1, v_2, v_3$, with unit axis $q$ and cyclic edge vectors $e_1, e_2, e_3$. Putting $x_i = v_i \cdot u$, the identity

$$\max_i x_i - \min_i x_i = 1/2 \sum_{\text{cyc}} |\, x_i - x_j|$$

gives

$$w_K(u) = 1/2 \sum_{k=1}^{3} |\, u \cdot e_k| + h\, |u \cdot q|$$

a sum of four absolute linear forms, although $K$ is not centrally symmetric and hence is not a zonotope.

Since $w_K = w_{DK}$, admitting bodies that are not centrally symmetric enlarges the class of bodies under consideration but not the class of width laws that arise; what the following records is the exact preimage of the zonotopal case under $K \mapsto DK$.

**Proposition 7.9.** Let $K$ be a full-dimensional convex body whose central symmetral $DK$ is the zonotope generated by pairwise non-parallel nonzero vectors $g_1, \dots, g_N$, and let $\mathcal{A}$ be the arrangement of the great circles $\{u \cdot g_k = 0\}$. On each open two-dimensional cell of $\mathcal{A}$ the function $w_K$ is linear. The law of $w_K$ is absolutely continuous, and for $\rho > 0$ its density is

$$f(\rho) = \frac{1}{4\pi} \sum_{\varepsilon} \frac{\Psi_\varepsilon(\rho)}{|G_\varepsilon|}, \qquad G_\varepsilon = \sum_{k=1}^{N} \varepsilon_k\, g_k,$$

where the sum is over the sign patterns $\varepsilon \in \{\pm 1\}^N$ whose cells are non-empty and two-dimensional. Here $\Psi_\varepsilon(\rho)$ is the angular measure of the part of the level circle $\{u \in S^2 : u \cdot G_\varepsilon = \rho\}$ lying in that cell, and is taken to be zero when the circle does not meet the cell, in particular when $\rho > |G_\varepsilon|$. The density is zero outside the range of $w_K$. Antipodal patterns $\varepsilon$ and $-\varepsilon$ give equal contributions. Moreover, the inclusion–exclusion calculation of $\Psi_\varepsilon$ terminates at order $N - 1$.

*Proof.* On the cell indexed by $\varepsilon$ one has $\varepsilon_k (u \cdot g_k) > 0$ for every $k$, and therefore

$$w_K(u) = \sum_{k=1}^{N} |u \cdot g_k| = u \cdot G_\varepsilon.$$

A non-empty cell cannot have $G_\varepsilon = 0$, since the displayed quantity is strictly positive there. For $0 < \rho < |G_\varepsilon|$, the level set $u \cdot G_\varepsilon = \rho$ is a circle of Euclidean radius

$$\frac{\sqrt{|G_\varepsilon|^2 - \rho^2}}{|G_\varepsilon|},$$

while the intrinsic gradient of $u \mapsto u \cdot G_\varepsilon$ on that circle has magnitude $\sqrt{|G_\varepsilon|^2 - \rho^2}$. The co-area quotient is consequently $1/|G_\varepsilon|$ times the angular measure of the part of the circle lying in the cell. Summing over the cells and dividing by the area $4\pi$ of $S^2$ gives the displayed density. The cell boundaries and endpoint level sets have spherical area zero, so they carry no mass, which proves absolute continuity.

To compute $\Psi_\varepsilon$, exclude from the full level circle the arcs on which at least one inequality $\varepsilon_k (u \cdot g_k) > 0$ fails. A point in the intersection of all $N$ excluded arcs would satisfy $\varepsilon_k (u \cdot g_k) \le 0$ for every $k$, and hence

$$\rho = u \cdot G_\varepsilon = \sum_{k=1}^{N} \varepsilon_k\, (u \cdot g_k) \le 0,$$

contrary to $\rho > 0$. The full $N$-fold intersection is therefore empty, so inclusion–exclusion terminates at order $N - 1$. Finally, the antipodal map sends the cell of $\varepsilon$ to that of $-\varepsilon$ and preserves the corresponding contribution. ▫

**Corollary 7.10.** Suppose $DK$ is the parallelepiped generated by $g_1, g_2, g_3$, so that $w_K = w_{DK}$, and write $V = |\det(g_1, g_2, g_3)|$; all diagonals, face areas and volumes below are those of $DK$. Then $N = 3$ and the inclusion-exclusion terminates at order two, as for a box. The eight cells form four antipodal pairs, and the density is a sum of four terms, one for each of the four space diagonals $\pm g_1 \pm g_2 \pm g_3$. The possible critical values of $w_K$ on the strata of $\mathcal{A}$ are

$$\frac{V}{|g_i \times g_j|} \quad \text{(0-strata)}, \qquad \left|P_k\left(\pm g_i \pm g_j\right)\right| \quad \text{(1-strata)}, \qquad |\pm g_1 \pm g_2 \pm g_3| \quad \text{(2-strata)}, \qquad \{i, j, k\} = \{1,2,3\},$$

with $P_k$ the orthogonal projection onto $g_k^{\perp}$: the volume divided by each face area, the projected face diagonals, and the four space diagonals. Their incidence is decided exactly as follows.

1. Each of the three 0-stratum values $V/\left|g_i \times g_j\right|$ is attained.
2. Fix $\{i, j, k\} = \{1,2,3\}$, put $p_i = P_k g_i$ and $p_j = P_k g_j$, and choose $\varepsilon_i, \varepsilon_j \in \{\pm 1\}$ modulo simultaneous reversal. The 1-stratum value

$$|\varepsilon_i p_i + \varepsilon_j p_j|$$

is attained in the relative interior of its sign arc if and only if

$$|p_i|^2 + \varepsilon_i \varepsilon_j\, p_i \cdot p_j > 0, \qquad \left|p_j\right|^2 + \varepsilon_i \varepsilon_j\, p_i \cdot p_j > 0.$$

Both conditions depend on the signs only through $\varepsilon_i \varepsilon_j$, so the two cases on this stratum are $\left|p_i + p_j\right|$ and $\left|p_i - p_j\right|$.

3. For $\varepsilon = (\varepsilon_1, \varepsilon_2, \varepsilon_3) \in \{\pm 1\}^3$ modulo simultaneous reversal, put $G_\varepsilon = \sum_r \varepsilon_r\, g_r$. The 2-stratum value $|G_\varepsilon|$ is attained in the relative interior of its sign cell if and only if

$$\varepsilon_r\, g_r \cdot G_\varepsilon > 0, \qquad r = 1,2,3.$$

In assertions 2 and 3, equality in one of the displayed conditions places the stationary direction on the boundary and transfers the value to the appropriate lower-dimensional stratum; a negative left side places it outside the closure of the specified sign cell. Thus these tests give the exact active list, with coincident numerical values combined only after the incidence test.

The density is continuous except possibly at endpoint values $|G_\varepsilon|$ for which the stationary direction belongs to the closure of the corresponding cell. To state the jump uniformly, suppose

$$\varepsilon_r\, g_r \cdot G_\varepsilon \geq 0, \qquad r = 1,2,3,$$

put $n = G_\varepsilon / |G_\varepsilon|$, and let $\Lambda_\varepsilon$ be the angular measure of the tangent directions $v \in S^1 \cap n^{\perp}$ satisfying

$$\varepsilon_r\, v \cdot g_r \geq 0 \quad \text{for every } r \text{ such that } \varepsilon_r\, g_r \cdot G_\varepsilon = 0.$$

The corresponding antipodal cell pair has a downward jump of size

$$\frac{\Lambda_\varepsilon}{2\pi|G_\varepsilon|}.$$

In the strict case of assertion 3, $\Lambda_\varepsilon = 2\pi$, giving the jump $1/|G_\varepsilon|$; when equality occurs, the partial step may superpose with the lower-stratum mechanism. Contributions add when endpoint values coincide. The support is

$$\operatorname{supp}(w_K) = \left[\frac{V}{\max_{i<j}|g_i \times g_j|},\ \max_\varepsilon |G_\varepsilon|\right].$$

*Proof.* The two points of the 0-stratum $g_i^\perp \cap g_j^\perp \cap S^2$ are $\pm(g_i \times g_j)/|g_i \times g_j|$. Since the generators are linearly independent, neither lies on $g_k^\perp$, and at either point

$$w_K(u) = |u \cdot g_k| = \frac{|\det(g_1, g_2, g_3)|}{|g_i \times g_j|}.$$

On the 1-stratum $g_k^\perp \cap S^2$, within the sign arc $\varepsilon_i(u \cdot g_i) > 0$, $\varepsilon_j(u \cdot g_j) > 0$, one has

$$w_K(u) = u \cdot (\varepsilon_i g_i + \varepsilon_j g_j) = u \cdot (\varepsilon_i p_i + \varepsilon_j p_j).$$

The stationary points of this linear form on the unit circle in $g_k^\perp$ are $\pm H/|H|$, where $H = \varepsilon_i p_i + \varepsilon_j p_j$; here $H \neq 0$ because $p_i, p_j$ are linearly independent. The minimum $-H/|H|$ cannot lie in the chosen arc, where $w_K > 0$. The stationary maximum $H/|H|$ belongs to its relative interior exactly when $\varepsilon_i p_i \cdot H > 0$ and $\varepsilon_j p_j \cdot H > 0$. Expanding these two inequalities gives assertion 2.

Finally, on the open 2-cell indexed by $\varepsilon$ one has $w_K(u) = u \cdot G_\varepsilon$. Its stationary maximum on $S^2$ is attained at $u = G_\varepsilon/|G_\varepsilon|$, and this direction belongs to the cell exactly when $\varepsilon_r(u \cdot g_r) > 0$ for every $r$. Multiplication by the positive number $|G_\varepsilon|$ gives the inequalities in assertion 3. Replacing a strict inequality by equality or making any left side negative gives the stated boundary and exterior alternatives.

The antipodal stationary minimum $-G_\varepsilon/|G_\varepsilon|$ cannot lie in the chosen cell, where $w_K > 0$. If the maximum is interior, then for $\rho < |G_\varepsilon|$ sufficiently close to $|G_\varepsilon|$ the entire shrinking level circle lies in the cell, and hence $\Psi_\varepsilon(\rho) \to 2\pi$. Proposition 7.9 then gives a limiting contribution $1/(2|G_\varepsilon|)$ from that cell and therefore $1/|G_\varepsilon|$ from its antipodal pair.

It remains to exclude other discontinuities. Fix a cell indexed by $\varepsilon$, write $G = G_\varepsilon$, $n = G/|G|$, and let $P$ denote orthogonal projection onto $n^\perp$. For $0 < \rho < |G|$, put $s = \sqrt{|G|^2 - \rho^2}$. On the level circle, the portion excluded by the $k$th sign inequality is an arc centred in the direction of $-\varepsilon_k P g_k$ and, when $P g_k \neq 0$, has half-width $\arccos c_k$, where

$$c_k = \frac{\rho\, \varepsilon_k (G \cdot g_k)}{|G|\, s\, |P g_k|}.$$

After clipping $c_k$ to $[-1,1]$, the centre and half-width depend continuously on $\rho$; the clipped values 1 and $-1$ represent the empty arc and the full circle. The angular length of a union of finitely many such arcs is continuous in their centres and half-widths. If $Pg_k = 0$, then $u \cdot g_k = (\rho/|G|)(n \cdot g_k)$ has constant sign on the level circle, so the corresponding constraint is either vacuous or fatal throughout and causes no discontinuity. Hence $\Psi_\varepsilon$ is continuous on $(0,|G|)$.

As $\rho \uparrow |G|$, the level circle shrinks to $n$. If some $\varepsilon_r(g_r \cdot G)$ is negative, the circle misses the cell near the endpoint and $\Psi_\varepsilon \to 0$. If all are non-negative, the strict inequalities become vacuous in the limit and every equality leaves the tangent half-circle condition $\varepsilon_r(v \cdot g_r) \geq 0$. Consequently $\Psi_\varepsilon \to \Lambda_\varepsilon$. Proposition 7.9 gives a left limit $\Lambda_\varepsilon/(4\pi|G|)$ from the cell and hence $\Lambda_\varepsilon/(2\pi|G|)$ from its antipodal pair; above $|G|$ that contribution is zero. Summation over the four pairs proves the asserted continuity and jump formula, and shows that there are no other discontinuities.

No minimum of $w_K$ lies in the relative interior of a 1- or 2-stratum: on either stratum the only admissible stationary point is the maximum just identified. The minimum is therefore attained on a 0-stratum, and assertion 1 gives the left endpoint of the support. Finally,

$$\max_{u \in S^2} w_K(u) = \max_{u \in S^2} \max_{\varepsilon} u \cdot G_\varepsilon = \max_{\varepsilon} |G_\varepsilon|,$$

which gives the right endpoint. ▫

For the orthogonal specialization $g_i = a_i e_i$, one has $p_i \cdot p_j = 0$, so every incidence inequality above is strict. The three lists reduce exactly to the edge lengths, face diagonals and space diagonal of Lemma 2.3.

**Proposition 7.11 (three-direction rigidity).** *If the central symmetral DK of a full-dimensional convex body $K \subset \mathbb{R}^3$ is a parallelepiped $Z$, then $K$ is a translate of $Z$.*

*Proof.* Write $Z = \sum_{i=1}^{3} [-1/2\, g_i, 1/2\, g_i]$, where the $g_i$ are linearly independent, and let $M$ be the matrix with columns $g_i$. The nonsingular map $M^{-1}$ commutes with central symmetrisation and sends $Z$ to $[-1/2, 1/2]^3$; it also preserves the assertion of being a translate. It therefore suffices to prove the result when $h_Z(u) = 1/2 \sum_i |u_i|$, which is additive on each closed orthant.

Put $p(u) = h_K(u) - h_Z(u)$. Since $h_K(u) + h_K(-u) = 2h_Z(u)$ and $h_Z$ is even, $p$ is odd. For $u, v$ in one orthant, the subadditivity of $h_K$ gives $p(u+v) \leq p(u) + p(v)$; applying the same inequality to $-u, -v$ and using oddness gives the reverse. Hence $p$ is additive and positively homogeneous on each orthant, so $p(u) = t_\varepsilon \cdot u$ there.

Adjacent orthants share a two-dimensional cone on which their two linear forms agree. Hence, if the orthants differ only in the sign of the $i$th coordinate,

$$t_\varepsilon - t_{\varepsilon'} = \lambda e_i$$

for some scalar $\lambda$, while oddness gives $t_{-\varepsilon} = t_\varepsilon$. Choose any three-step path in the orthant graph from $\varepsilon$ to $-\varepsilon$, flipping each coordinate once. Summing the three transition relations along this path gives

$$0 = t_{-\varepsilon} - t_\varepsilon = \lambda_1 e_1 + \lambda_2 e_2 + \lambda_3 e_3,$$

and the independence of $e_1, e_2, e_3$ forces $\lambda_1 = \lambda_2 = \lambda_3 = 0$. Since any prescribed adjacency edge can be included in such a three-step path, every transition coefficient vanishes. Thus the linear forms agree across all orthants, $p(u) = t \cdot u$ globally, and therefore $h_K = h_Z + t \cdot u$, which is equivalent to $K = Z + t$. ▫

The enlargement afforded by Proposition 7.7 first becomes genuine when an irredundant representation of $DK$ involves at least four distinct generator directions, as in Remark 7.8; the number of generators alone is not invariant, since one may always be split into two parallel ones.

**Remark 7.12.** Orthogonality makes the box exceptional in two ways: the sign-change symmetries act transitively on its cells, so a single term suffices, and its four space-diagonal vectors have equal length, so three of the four cell endpoints are hidden at the upper endpoint of the support. For a non-orthogonal parallelepiped these symmetries generally fail and the four candidate lengths may differ. Corollary 7.10 decides which stationary directions lie in the cell closures and gives the resulting full or partial steps exactly. The corresponding classification of the corner and fold mechanisms for a general zonotope is Conjecture 7.13.

**Conjecture 7.13.** Let $DK$ be a zonotope with generators $g_1, \dots, g_N$, reduced as follows: within each parallel class choose a common orientation, replacing $g_k$ by $-g_k$ where necessary, which changes neither $|u \cdot g_k|$ nor the segment $[-1/2\, g_k, 1/2\, g_k]$, and then replace the class by the single generator equal to the sum of its members. Let $\mathcal{A}$ be the arrangement of the great circles $\{u \cdot g_k = 0\}$ of the reduced family. The reduction is needed because $DK$ and $w_K$ are unchanged by it, whereas an unreduced family would make $\mathcal{A}$ depend on a redundant representation. Then the density of $w_K$ is real-analytic away from the critical values of the restrictions of $w_K$ to the strata of $\mathcal{A}$ and the endpoint values of the two-dimensional cells. A 0-stratum contributes a corner, a stationary value on a 1-stratum contributes a square-root fold, and an endpoint of a 2-cell contributes a step whenever its stationary direction lies in the closure of that cell. If several mechanisms occur at the same value, they superpose, the strongest being dominant in the order

$$\text{step discontinuity} \succ \text{square-root fold} \succ \text{corner}$$

Theorem 3.8 is the case of a box with $a_1 = d_1$, where a fold and a corner superpose. The local corner and fold mechanisms underlying Theorems 3.4 and 3.7 do not depend on orthogonality; what is required in general is the determination, near each critical value, of which excluded arcs are non-empty and which pairs of them overlap, supplied for the box by Proposition 3.3, together with the incidence condition of Corollary 7.10.

# 8. Concluding remarks

Walters's projected-area calculation establishes the underlying rectangular-parallelepiped law. The present treatment supplies the global positive-part and CDF forms in width coordinates and then develops what those forms reveal: the complete singularity structure, the cumulant degree and closure theorems, inverse reconstruction and the boundary between box-specific and general convex-body behaviour.

We record three directions of extension, the third of which accounts for the restriction to three dimensions.

**Zonotopes.** Corollary 7.10 completes the incidence analysis for three independent generator directions. Conjecture 7.13 asserts that the classification of Theorem B extends to every body whose central symmetral is a zonotope; what remains is to extend the explicit incidence criteria to arbitrary reduced generator families and prove the local mechanism at each attained value.

**General polytopes.** For an arbitrary full-dimensional polytope $P$, the width $w_P$ is piecewise linear on $S^2$, being linear on each cell of the common refinement of the normal fans of $P$ and $-P$. The argument of Theorem 3.1 therefore applies cell by cell and gives an absolutely continuous law whose density is a finite sum of angular measures. What is lost is the global indexing by a single family of generator sign patterns, and with it the uniform generator-based formula of Proposition 7.9. Each polyhedral cell is nevertheless cut out by finitely many linear inequalities and therefore admits a finite calculation, although the resulting indexing may be considerably less economical. Whether a usable closed form survives is Question 8.3.

**Higher dimensions.** The cancellation that produces Theorem 3.1 is a coincidence of dimension three. For $u$ uniform on $S^{n-1}$ and $w = \sum_{i=1}^{n} a_i\,|u_i|$, the function $w$ is linear on each open orthant; there the level set $\{w = \rho\}$ is the portion, cut out by that orthant, of an $(n-2)$-sphere of radius $s/d$, where $s = \sqrt{d^2 - \rho^2}$ and $d = |a|$, while the intrinsic gradient of $w$ has magnitude $s$. The co-area formula therefore contributes the factor

$$\frac{1}{s}\Big(\frac{s}{d}\Big)^{n-2} = \frac{s^{\,n-3}}{d^{\,n-2}}$$

which is independent of $\rho$ precisely when $n = 3$. In the plane the surviving factor $s^{-1}$ is visible in Proposition 3.12. For $n > 3$ a factor $s^{\,n-3}$ remains, and the angular measure becomes the volume of a region of an $(n-2)$-sphere cut by $n$ hyperplanes; the argument of Proposition 7.9 still terminates the inclusion-exclusion at order $n-1$, but the individual terms are no longer arccosines.

The two halves of the paper combine into a single statement of what the law of the width does and does not retain.

**Proposition 8.1.** *Within the class of rectangular boxes, the law of $w_K$ determines $K$ up to congruence: $\kappa_1$ gives $L$, $\kappa_2$ gives $S$, $\kappa_3$ gives $V$, and Lemma 2.1 returns the edge lengths up to permutation. Within the class of all convex bodies the law does not determine $K$ up to congruence. More strongly, if $K$ is not centrally symmetric, then $K$ and its central symmetral $DK = 1/2\,(K - K)$ are non-congruent bodies with the same width function pointwise.*

*Proof.* The assertion for boxes is Theorem 6.8: the first three cumulants determine $L, S, V$, and Lemma 2.1 then recovers the unordered triple of edge lengths.

For an arbitrary convex body,

$$h_{DK}(u) = 1/2\,h_{K-K}(u) = 1/2\,\big(h_K(u) + h_K(-u)\big) = 1/2\,w_K(u)$$

Since $DK$ is centrally symmetric, its support function is even, and therefore

$$w_{DK}(u) = h_{DK}(u) + h_{DK}(-u) = 2h_{DK}(u) = w_K(u)$$

for every $u \in S^2$. If $K$ is not centrally symmetric, then it cannot be congruent to $DK$, because central symmetry is preserved by every Euclidean isometry whereas $DK$ is centrally symmetric. Thus even the pointwise width function, and hence certainly its law, does not determine a general convex body up to congruence. ▫

Bodies of a fixed constant width provide a particularly large family with the same degenerate width law, although this is another manifestation of the common central symmetral rather than a failure within the centrally symmetric class: a centrally symmetric body of constant width $c$ has $h_K$ even, so $c = w_K = 2h_K$ and $K$ is the ball of radius $c/2$.

We record several questions that these considerations leave open.

**Question 8.2.** Between these two extremes, in which classes of bodies does the law of $w_K$ determine $K - K$ up to orthogonal equivalence? The law retains less than the pointwise width function, discarding the directional labelling of $h_{K-K}$ and being invariant under orthogonal transformations. For comparison, Aleksandrov's projection theorem [2] states that two centrally symmetric convex bodies with non-empty interior in $\mathbb{R}^n$ whose projections onto every hyperplane have equal volume are translates of one another; that concerns pointwise projection data, whereas the law of $w_K$ retains only the distribution of the values of a support function.

**Question 8.3.** For a general convex polytope $P$, does the cellwise inverse-trigonometric calculation of the density of $w_P$ admit a compact global representation, or one whose algebraic and combinatorial complexity can be bounded effectively in terms of the face data of $P$? The second paragraph above gives a finite cellwise procedure, but the economical sign-pattern indexing of Proposition 7.9 generally does not survive.

**Question 8.4.** What replaces the arccosines in dimensions above three, where the co-area factor $s^{n-3}$ no longer cancels?

**Question 8.5.** The third moment of the brightness of a general convex body is available by Corollary 7.3, because the trivariate Gaussian absolute product moment has a finite expression in algebraic and inverse-trigonometric functions [10], [8]. For four variables Nabeya [11] gives $\mathbb{E}|X_1X_2X_3X_4|$ as a finite combination of such functions together with the residual quantity $S(R) = \mathbb{E}[\mathrm{sgn}(X_1X_2X_3X_4)]$, which is also the non-pairwise term in the quadrivariate orthant probability and is represented there by a one-dimensional integral rather than in finite terms. In the case needed here the $X_i = g \cdot v_i$ arise from unit vectors $v_i \in \mathbb{R}^3$, so $R$ is a Gram matrix of rank at most three and $\det R = 0$. The sign of the product is then constant on each cell cut out on $S^2$ by the four great circles $v_i^{\perp}$. After the incidence type has been fixed, $S(R)$ is therefore a signed sum of the areas of finitely many spherical polygons, each of which has a finite expression in terms of its vertex angles. What is not clear is whether these casewise expressions can be combined into a compact formula in the six correlations without enumerating the incidence types. Is there such a formula? An affirmative answer would give the fourth brightness moment for every convex body in the same form as Theorem 7.1 and Proposition 7.2. Nabeya's derivation assumes $R$ non-singular, but the extension to the Gram case is immediate: for $X^{(\varepsilon)} = \sqrt{1-\varepsilon}\,X + \sqrt{\varepsilon}\,Z$ with $Z$ standard normal in $\mathbb{R}^4$ independent of $X$, the correlation matrix $R_\varepsilon = (1-\varepsilon)R + \varepsilon I_4$ is non-singular for $\varepsilon > 0$, and $X^{(\varepsilon)} \to X$ almost surely. The absolute product moments converge by uniform integrability, and so does the sign term: each $X_i$ is a non-degenerate Gaussian, so $\Pr[X_i = 0] = 0$, hence $\mathrm{sgn}\left(\prod_i X_i^{(\varepsilon)}\right) \to \mathrm{sgn}(\prod_i X_i)$ almost surely, and dominated convergence gives $S(R_\varepsilon) \to S(R)$. Nabeya's identity may therefore be passed to the singular limit along this regularisation. In degenerate subcases, individual partial-correlation terms should be interpreted through the same limiting procedure rather than assumed to be separately defined at $\det R = 0$.

**Question 8.6.** Among boxes subject to a prescribed constraint, such as fixed volume, surface area, mean width or diameter, which extremise a given functional of the width law, for instance its variance, its coefficient of variation, its entropy or its upper tail mass? Corollary 5.9 answers one instance, the cube being the unique minimiser of the coefficient of variation. Of the functionals named, the variance and the coefficient of variation are defined without further convention; an entropy or a tail functional requires one, for instance the differential entropy of $w/\mathbb{E}[w]$ or the mass $\Pr\left[w \geq \lambda\, \mathbb{E}[w]\right]$ at a fixed $\lambda > 1$. Once the normalisation and $\lambda$ are fixed, Proposition 3.11 evaluates this tail exactly.

In the planar setting, Akiyama and Kamae [1] determine the extremisers of the corresponding width-deviation ratio for convex $n$-gons. Their result provides a natural comparison for the three-dimensional optimisation problem posed here, but does not address rectangular boxes in $\mathbb{R}^3$.

**Question 8.7.** Is there a treatment covering the box and the ellipsoid uniformly? The box and the ellipsoid are two explicitly tractable $\Pi$-stable subfamilies, their widths being respectively a linear form in $|u|$ and the square root of a quadratic form in $u$; the projection bodies are precisely the centred zonoids [5, §4.1]. For the box the mean width is $L/2$ and the density is Theorem 3.1; for the ellipsoid the mean width is an elliptic integral and the density of the quadratic form is the object treated by Hillier. The two densities are piecewise on intervals separated respectively by the edge lengths and face diagonals in the box case, and by the characteristic roots in the ellipsoid case.

We close with an observation that motivates the sequel to this paper. It is a bridge rather than a result of the present one: it uses only the upper endpoint of the support, and none of the density, moment or cumulant theory developed above.

**A spectral family of width laws.** For $a_1, a_2, a_3 > 0$ and a nonzero index $n = (n_1, n_2, n_3)$ with $n_i \geq 0$, put

$$b^{(n)}(a) = \left(\frac{n_1}{a_1}, \frac{n_2}{a_2}, \frac{n_3}{a_3}\right)$$

and define

$$W_{n,a}(u) = \sum_{i=1}^{3} \frac{n_i}{a_i} |u_i|, \qquad u \in S^2$$

By Cauchy-Schwarz,

$$\max_{u \in S^2} W_{n,a}(u) = \Big(\sum_{i=1}^{3} \frac{n_i^2}{a_i^2}\Big)^{1/2}$$

Consequently, the Dirichlet and Neumann eigenvalues of the box satisfy

$$\lambda_n^D(a) = \pi^2 (\max W_{n,a})^2, \qquad n \in \mathbb{Z}_{>0}^3$$

and

$$\lambda_n^N(a) = \pi^2(\max W_{n,a})^2, \qquad n \in \mathbb{Z}_{\geq 0}^3, \quad n \neq (0,0,0)$$

respectively. When every $n_i > 0$, $W_{n,a}$ is the width of the auxiliary box with edge lengths $n_i/a_i$, so the density, moment and cumulant formulas of this paper apply directly. If one or more coordinates vanish, the corresponding laws are obtained as degenerate limits and describe auxiliary rectangles or segments.

Thus a fixed box gives the indexed family of auxiliary width laws

$$\{\mathcal{L}(W_{n,a})\}_n.$$

For each separated mode, the corresponding eigenvalue is $\pi^2$ times the square of the upper endpoint of the associated law. The Dirichlet or Neumann spectrum itself is the multiset of these squared upper endpoints, with multiplicities; it retains neither the remaining distributional information nor, in general, the mode labels. The possible use of that additional distributional information is a separate problem.

---

# Appendix A. The moment and cumulant algebra

This appendix carries out the algebra behind Propositions 4.8 and 5.6. Every coefficient through order five is obtained from the finite moment sum (6), and the cumulants are then collected explicitly by powers of $\pi^{-1}$.

**A.1 The spherical moments.** Corollary 4.2 gives, for $Q = \sum_i q_i$,

$$\mathbb{E}\prod_{i=1}^{3}|u_i|^{q_i} = \frac{1}{2\pi\,\Gamma(3+Q/2)}\prod_{i=1}^{3}\Gamma\left(\frac{q_i+1}{2}\right)$$

Only

$$\Gamma(1/2) = \sqrt{\pi}, \quad \Gamma(1) = 1, \quad \Gamma(3/2) = \sqrt{\pi}/2\,, \quad \Gamma(2) = 1,$$

$$\Gamma(5/2) = 3\sqrt{\pi}/4\,, \quad \Gamma(3) = 2, \quad \Gamma(7/2) = 15\sqrt{\pi}/8\,, \quad \Gamma(4) = 6$$

are required through order five. A factor $\sqrt{\pi}$ appears in $\Gamma\big((q+1)/2\big)$ exactly when $q$ is even, including $q = 0$, and in $\Gamma\big((3+Q)/2\big)$ exactly when $Q$ is even. If $k$ of the three exponents $q_i$ are even, then the numerator contributes $\pi^{k/2}$ and the denominator contributes

$$\pi^{1+\frac{1}{2}[Q \text{ even}]}$$

Since $Q \equiv 3 - k \pmod 2$, the cases $k = 3,2,1,0$ give respectively $\pi^0, \pi^0, \pi^{-1}, \pi^{-1}$. Thus a term carries $\pi^{-1}$ exactly when two or three of the $q_i$ are odd, and carries $\pi^0$ when zero or one is odd. This is the parity count of Lemma 4.5.

**A.2 The multinomial expansion.** Expanding

$$w^N = (\sum_i a_i\,|u_i|)^N$$

and applying A.1 termwise gives (6). Group the exponent triples into types, a type being an unordered partition $q$ of $N$ into three non-negative parts. Let $\sigma_q$ be the sum of the monomials $\prod_i a_i^{q_i}$ over the distinct assignments of that type. The coefficient of $\sigma_q$ is

$$c_q = \frac{N!}{2\pi\,\Gamma(N+3/2)} \prod_{i=1}^{3} \frac{\Gamma(q_i+1/2)}{q_i!}$$

so that

$$m_N = \sum_q c_q\,\sigma_q$$

The complete list through $N = 5$ is as follows.

| $N$ | **type** $q$ | **odd parts** | $c_q$ | $\sigma_q$ **in** $e_1, e_2, e_3$ |
|---|---|---|---|---|
| 1 | (1,0,0) | 1 | $1/2$ | $e_1$ |
| 2 | (2,0,0) | 0 | $1/3$ | $e_1^2 - 2e_2$ |
| 2 | (1,1,0) | 2 | $4/3\pi$ | $e_2$ |
| 3 | (3,0,0) | 1 | $1/4$ | $e_1^3 - 3e_1e_2 + 3e_3$ |
| 3 | (2,1,0) | 1 | $3/8$ | $e_1e_2 - 3e_3$ |
| 3 | (1,1,1) | 3 | $3/2\pi$ | $e_3$ |
| 4 | (4,0,0) | 0 | $1/5$ | $e_1^4 - 4e_1^2e_2 + 2e_2^2 + 4e_1e_3$ |
| 4 | (3,1,0) | 2 | $16/15\pi$ | $e_1^2e_2 - 2e_2^2 - e_1e_3$ |
| 4 | (2,2,0) | 0 | $2/5$ | $e_2^2 - 2e_1e_3$ |
| 4 | (2,1,1) | 2 | $8/5\pi$ | $e_1e_3$ |
| 5 | (5,0,0) | 1 | $1/6$ | $e_1^5 - 5e_1^3e_2 + 5e_1e_2^2 + 5e_1^2e_3 - 5e_2e_3$ |
| 5 | (4,1,0) | 1 | $5/16$ | $e_1^3e_2 - 3e_1e_2^2 - e_1^2e_3 + 5e_2e_3$ |
| 5 | (3,2,0) | 1 | $5/12$ | $e_1e_2^2 - 2e_1^2e_3 - e_2e_3$ |
| 5 | (3,1,1) | 3 | $5/3\pi$ | $e_3(e_1^2 - 2e_2)$ |
| 5 | (2,2,1) | 1 | $5/8$ | $e_2e_3$ |

Multiplying $c_q$ by $\sigma_q$ and summing over the rows of each order gives

$$m_1 = \frac{e_1}{2}$$

$$m_2 = \frac{e_1^2 - 2e_2}{3} + \frac{4e_2}{3\pi}$$

and

$$\begin{aligned} m_3 \quad &= 1/4\,(e_1^3 - 3e_1e_2 + 3e_3) + 3/8\,(e_1e_2 - 3e_3) + \frac{3e_3}{2\pi} \\ &= \frac{2e_1^3 - 3e_1e_2 - 3e_3}{8} + \frac{3e_3}{2\pi} \end{aligned}$$

At order four,

$$\begin{aligned} m_4 = \quad & 1/5\,(e_1^4 - 4e_1^2e_2 + 2e_2^2 + 4e_1e_3) + 2/5\,(e_2^2 - 2e_1e_3) \\ & + \frac{16}{15\pi}(e_1^2e_2 - 2e_2^2 - e_1e_3) + \frac{8e_1e_3}{5\pi} \end{aligned}$$

The rational terms in $e_1e_3$ cancel, and the remaining terms collect to

$$m_4 = \frac{e_1^4 - 4e_1^2e_2 + 4e_2^2}{5} + \frac{8(2e_1^2e_2 - 4e_2^2 + e_1e_3)}{15\pi}$$

At order five,

$$\begin{aligned} m_5 = \quad & 1/6\left(e_1^5 - 5e_1^3e_2 + 5e_1e_2^2 + 5e_1^2e_3 - 5e_2e_3\right) \\ & + 5/16\,(e_1^3e_2 - 3e_1e_2^2 - e_1^2e_3 + 5e_2e_3) \\ & + 5/12\,(e_1e_2^2 - 2e_1^2e_3 - e_2e_3) + 5/8\,e_2e_3 + \frac{5e_3(e_1^2 - 2e_2)}{3\pi} \end{aligned}$$

Putting the rational terms over the common denominator 48 gives the coefficient vector

$$\frac{1}{48}(8, -25, 15, -15, 45)$$

in the ordered basis

$$e_1^5, \quad e_1^3e_2, \quad e_1e_2^2, \quad e_1^2e_3, \quad e_2e_3$$

and therefore

$$m_5 = \frac{8e_1^5 - 25e_1^3e_2 + 15e_1e_2^2 - 15e_1^2e_3 + 45e_2e_3}{48} + \frac{5e_3(e_1^2 - 2e_2)}{3\pi}$$

These are precisely the formulas of Proposition 4.8.

**A.3 The symmetric reductions.** Write $p_j = \sum_i a_i^j$. Newton's identities give

$$p_1 = e_1, \qquad p_2 = e_1^2 - 2e_2, \qquad p_3 = e_1^3 - 3e_1e_2 + 3e_3$$

$$p_4 = e_1^4 - 4e_1^2e_2 + 2e_2^2 + 4e_1e_3$$

$$p_5 = e_1^5 - 5e_1^3e_2 + 5e_1e_2^2 + 5e_1^2e_3 - 5e_2e_3$$

The remaining reductions used in A.2 are

$$\sigma_{(1,1,0)} = e_2$$

$$\sigma_{(N-1,1,0)} = p_{N-1}p_1 - p_N \qquad (N \geq 3)$$

$$\sigma_{(2,2,0)} = 1/2\,(p_2^2 - p_4), \qquad \sigma_{(3,2,0)} = p_3p_2 - p_5$$

$$\sigma_{(2,1,1)} = e_3e_1, \qquad \sigma_{(3,1,1)} = e_3p_2, \qquad \sigma_{(2,2,1)} = e_3e_2$$

For example,

$$\sigma_{(3,1,0)} = \sum_{i\neq j} a_i^3\, a_j = \left(\sum_i a_i^3\right)\left(\sum_j a_j\right) - \sum_i a_i^4 = p_3p_1 - p_4$$

The restriction $N \geq 3$ in the second identity is essential: for $N = 2$, the two exponents coincide and $p_1^2 - p_2 = 2e_2$, whereas $\sigma_{(1,1,0)} = e_2$.

**A.4 The cumulants.** Put

$$A = m_1 = \frac{e_1}{2}$$

and decompose each higher moment as

$$m_j = R_j + \frac{P_j}{\pi}$$

where

$$R_2 = \frac{e_1^2 - 2e_2}{3}, \qquad P_2 = \frac{4e_2}{3}$$

$$R_3 = \frac{2e_1^3 - 3e_1e_2 - 3e_3}{8}, \qquad P_3 = \frac{3e_3}{2}$$

$$R_4 = \frac{e_1^4 - 4e_1^2e_2 + 4e_2^2}{5}, \qquad P_4 = \frac{8(2e_1^2e_2 - 4e_2^2 + e_1e_3)}{15}$$

$$R_5 = \frac{8e_1^5 - 25e_1^3e_2 + 15e_1e_2^2 - 15e_1^2e_3 + 45e_2e_3}{48}, \qquad P_5 = \frac{5e_3(e_1^2 - 2e_2)}{3}$$

The standard moment-cumulant identities are

$$\kappa_2 = m_2 - m_1^2$$

$$\kappa_3 = m_3 - 3m_2m_1 + 2m_1^3$$

$$\kappa_4 = m_4 - 4m_3m_1 - 3m_2^2 + 12m_2m_1^2 - 6m_1^4$$

$$\kappa_5 = m_5 - 5m_4m_1 - 10m_3m_2 + 20m_3m_1^2 + 30m_2^2m_1 - 60m_2m_1^3 + 24m_1^5$$

For the second cumulant,

$$\begin{aligned}\kappa_2 &= R_2 - A^2 + \frac{P_2}{\pi} \\ &= \frac{e_1^2}{12} - \frac{2e_2}{3} + \frac{4e_2}{3\pi}\end{aligned}$$

For the third,

$$\kappa_3 = (R_3 - 3AR_2 + 2A^3) + \frac{P_3 - 3AP_2}{\pi}$$

and the two coefficients are

$$R_3 - 3AR_2 + 2A^3 = \frac{5e_1e_2}{8} - \frac{3e_3}{8}$$

$$P_3 - 3AP_2 = -2e_1e_2 + \frac{3e_3}{2}$$

Thus

$$\kappa_3 = \frac{5e_1e_2}{8} - \frac{2e_1e_2}{\pi} - \frac{3e_3}{8} + \frac{3e_3}{2\pi}$$

For the fourth cumulant, separate the three powers of $\pi^{-1}$:

$$\kappa_4 = K_{4,0} + \frac{K_{4,1}}{\pi} + \frac{K_{4,2}}{\pi^2}$$

where

$$K_{4,0} = R_4 - 4AR_3 - 3R_2^2 + 12A^2R_2 - 6A^4$$

$$K_{4,1} = P_4 - 4AP_3 - 6R_2P_2 + 12A^2P_2$$

$$K_{4,2} = -3P_2^2$$

Substitution gives

$$\begin{aligned}K_{4,0} = &\quad (1/5 - 1/2 - 1/3 + 1 - 3/8)e_1^4 \\ &\quad +(-4/5 + 3/4 + 4/3 - 2)e_1^2e_2 + (4/5 - 4/3)e_2^2 + 3/4\, e_1e_3 \\ = &\quad -\frac{e_1^4}{120} - \frac{43e_1^2e_2}{60} - \frac{8e_2^2}{15} + \frac{3e_1e_3}{4}\end{aligned}$$

$$\begin{aligned}K_{4,1} = &\quad (16/15 - 8/3 + 4)e_1^2e_2 + (-32/15 + 16/3)e_2^2 + (8/15 - 3)e_1e_3 \\ = &\quad \frac{12e_1^2e_2}{5} + \frac{16e_2^2}{5} - \frac{37e_1e_3}{15}\end{aligned}$$

and

$$K_{4,2} = -3(\frac{4e_2}{3})^2 = -\frac{16e_2^2}{3}$$

Therefore

$$\kappa_4 = -\frac{e_1^4}{120} - \frac{43e_1^2 e_2}{60} + \frac{12e_1^2 e_2}{5\pi} + \frac{3e_1 e_3}{4} - \frac{37e_1 e_3}{15\pi} - \frac{8e_2^2}{15} + \frac{16e_2^2}{5\pi} - \frac{16e_2^2}{3\pi^2}$$

For the fifth cumulant, write

$$\kappa_5 = K_{5,0} + \frac{K_{5,1}}{\pi} + \frac{K_{5,2}}{\pi^2}$$

The rational part is

$$K_{5,0} = R_5 - 5AR_4 - 10R_3R_2 + 20A^2R_3 + 30AR_2^2 - 60A^3R_2 + 24A^5$$

Its coefficient collection is shown in the ordered basis

$$e_1^5, \quad e_1^3 e_2, \quad e_1 e_2^2, \quad e_1^2 e_3, \quad e_2 e_3$$

| **contribution** | $e_1^5$ | $e_1^3 e_2$ | $e_1 e_2^2$ | $e_1^2 e_3$ | $e_2 e_3$ |
|---|---|---|---|---|---|
| $R_5$ | $1/6$ | $-25/48$ | $5/16$ | $-5/16$ | $15/16$ |
| $-5AR_4$ | $-1/2$ | $2$ | $-2$ | $0$ | $0$ |
| $-10R_3R_2$ | $-5/6$ | $35/12$ | $-5/2$ | $5/4$ | $-5/2$ |
| $20A^2R_3$ | $5/4$ | $-15/8$ | $0$ | $-15/8$ | $0$ |
| $30AR_2^2$ | $5/3$ | $-20/3$ | $20/3$ | $0$ | $0$ |
| $-60A^3R_2$ | $-5/2$ | $5$ | $0$ | $0$ | $0$ |
| $24A^5$ | $3/4$ | $0$ | $0$ | $0$ | $0$ |
| total | $0$ | $41/48$ | $119/48$ | $-15/16$ | $-25/16$ |

Hence

$$K_{5,0} = \frac{41e_1^3 e_2}{48} + \frac{119e_1 e_2^2}{48} - \frac{15e_1^2 e_3}{16} - \frac{25e_2 e_3}{16}$$

The coefficient of $\pi^{-1}$ is

$$K_{5,1} = P_5 - 5AP_4 - 10(R_3P_2 + P_3R_2) + 20A^2P_3 + 60AR_2P_2 - 60A^3P_2$$

Its collection in the ordered basis

$$e_1^3 e_2, \quad e_1 e_2^2, \quad e_1^2 e_3, \quad e_2 e_3$$

is

| **contribution** | $e_1^3 e_2$ | $e_1 e_2^2$ | $e_1^2 e_3$ | $e_2 e_3$ |
|---|---|---|---|---|
| $P_5$ | 0 | 0 | $5/3$ | $-10/3$ |
| $-5AP_4$ | $-8/3$ | $16/3$ | $-4/3$ | 0 |
| $-10(R_3P_2 + P_3R_2)$ | $-10/3$ | 5 | $-5$ | 15 |
| $20A^2P_3$ | 0 | 0 | $15/2$ | 0 |
| $60AR_2P_2$ | $40/3$ | $-80/3$ | 0 | 0 |
| $-60A^3P_2$ | $-10$ | 0 | 0 | 0 |
| total | $-8/3$ | $-49/3$ | $17/6$ | $35/3$ |

Therefore

$$K_{5,1} = -\frac{8e_1^3e_2}{3} - \frac{49e_1e_2^2}{3} + \frac{17e_1^2e_3}{6} + \frac{35e_2e_3}{3}$$

Finally,

$$\begin{aligned} K_{5,2} &= -10P_3P_2 + 30AP_2^2 \\ &= -10\left(\frac{3e_3}{2}\right)\left(\frac{4e_2}{3}\right) + 30\left(\frac{e_1}{2}\right)(\frac{4e_2}{3})^2 \\ &= -20e_2e_3 + \frac{80e_1e_2^2}{3} \end{aligned}$$

Combining the three powers gives

$$\begin{aligned} \kappa_5 = \quad & \frac{41e_1^3e_2}{48} - \frac{8e_1^3e_2}{3\pi} - \frac{15e_1^2e_3}{16} + \frac{17e_1^2e_3}{6\pi} \\ & + \frac{119e_1e_2^2}{48} - \frac{49e_1e_2^2}{3\pi} + \frac{80e_1e_2^2}{3\pi^2} \\ & - \frac{25e_2e_3}{16} + \frac{35e_2e_3}{3\pi} - \frac{20e_2e_3}{\pi^2} \end{aligned}$$

This completes the hand derivation of the cumulants in Proposition 5.6.

**A.5 Mean normalisation.** Set $e_1 = 1$. Then

$$e_2 = \frac{1-G}{2}, \qquad e_3 = J, \qquad \kappa_1 = \frac{1}{2}$$

and therefore

$$g_n = \frac{\kappa_n}{\kappa_1^n} = 2^n \kappa_n$$

Substitution into the formulas of A.4 gives

$$g_2 = \frac{4G-3}{3} + \frac{8(1-G)}{3\pi}$$

$$g_3 = -\frac{5G+6J-5}{2} + \frac{4(2G+3J-2)}{\pi}$$

$$g_4 = -\frac{2(16G^2-75G-90J+60)}{15} + \frac{16(12G^2-42G-37J+30)}{15\pi} - \frac{64(1-G)^2}{3\pi^2}$$

$$g_5 = \frac{119G^2+150GJ-320G-330J+201}{6} - \frac{8(49G^2+70GJ-114G-104J+65)}{3\pi} + \frac{320(G-1)(2G+3J-2)}{3\pi^2}$$

which are the mean-normalised formulas of Proposition 5.6. The factor $2^n$ is essential: setting $e_1 = 1$ normalises the edge sum, whereas $g_n$ is normalised by the mean $\kappa_1 = e_1/2$.

**A.6 A worked check.** For the unit cube, $e_1 = 3$, $e_2 = 3$, $e_3 = 1$. Hence

$$m_1 = \frac{3}{2}$$

and

$$m_2 = \frac{9-6}{3} + \frac{12}{3\pi} = 1 + \frac{4}{\pi}$$

so

$$\kappa_2 = m_2 - m_1^2 = \frac{4}{\pi} - \frac{5}{4}$$

Remark 4.11 reproduces $m_1$ and $m_2$ by numerical integration of the independently derived density of Theorem 3.1; that computation confirms the moment formulas and is not used in their derivation. The density normalisation itself is exact by Proposition 3.11.

---

**Acknowledgements**. The material in this paper has been in development for many years and appears in print only now. I thank my family for their patience over that time. I am grateful to the early readers who criticised successive drafts; in particular, an error in an earlier version of the argument of Section 6 was caught by a careful reader. I owe a lasting debt to my late mathematics teacher, Amos Matalon, who showed me advanced mathematics at an early age, insisted that it be taken seriously, and set me on the path to research. Any remaining errors are my own.

**Use of AI tools**. All derivations in this paper were carried out by hand. A large language model was used as an auxiliary tool in two respects: to run independent numerical checks confirming the closed forms stated here, and to assist with prior-art and bibliographic searching. Figures 1 and 2 are deterministic plots, produced with matplotlib from the formulas proved in this paper; the plotting script was drafted with the assistance of a large language model and checked by the author. No figure in this paper is AI-generated imagery. The model produced no text and no mathematical argument in the paper, and the author takes full responsibility for its content.

**Declarations**. Funding: this research received no external funding. Competing interests: the author declares none. Data availability: the script generating Figures 1 and 2 is archived at [10.5281/zenodo.21780025]. It is used for plotting and numerical confirmation only; no formula in this paper is derived from it.

---